\documentclass{amsart}

\usepackage[usenames,dvipsnames,svgnames,table]{xcolor}

\usepackage{fancyvrb}

\usepackage{geometry}      
\usepackage[normalem]{ulem}          
\usepackage{wrapfig}
\usepackage{graphics}
\usepackage{graphicx}
\usepackage{mathrsfs}
\usepackage{amssymb}
\usepackage{amsfonts}
\usepackage{mathrsfs}
\usepackage{epstopdf}
\usepackage{lscape}
\usepackage{multicol}
\usepackage[utf8]{inputenc}
\usepackage{tikz,caption}
\usepackage{enumitem}
\setlist[itemize]{leftmargin=2em}
\setlist[enumerate]{leftmargin=2em}
\usepackage{booktabs}
\usepackage{multirow}
\usepackage{mathtools}
\usepackage[linesnumbered,ruled]{algorithm2e}

\usepackage{soul}
\usepackage{upgreek}

\title{
A memorial tribute
\\ Adriano Garsia (1928--2024)
}

\author[Barcelo, et al.]{H\'el\`ene Barcelo, Fran\c{c}ois Bergeron,
  Nantel Bergeron, Sara Billey,
  Anders Bj\"orner, Tullio Ceccherini-Silberstein,
  Michele D'Adderio, Mark Haiman, Angela Hicks, Luc Lapointe,
  David Little,
  Jennifer Morse, Arun Ram, Marino Romero, Emily Sergel,
  Michelle Wachs, Guoce Xin, Mike Zabrocki
}

\date{}

\begin{document}

\maketitle

\section*{Introduction}

Adriano Mario Garsia was born in Tunis on August 20, 1928,
to a Tunisian-Italian family. He lived on a farm there until
the end of World War II, then moved to Rome.  After finishing 
high school, he was sent to the United States to live with relatives
in Woyming
and eventually made his way to California, becoming a student of
Charles Loewner at Stanford in the early 1950s. Following his Ph.D.,
Adriano held positions at MIT, the University of Minnesota, and Caltech
before joining the nascent mathematics department at 
the University of California, San Diego, in 1966 where he spent the remainder of his career. 
He passed away in San Diego on October 6, 2024, at the age of 96.

  \begin{wrapfigure}{r}{0.4\textwidth}
    \centering
    \includegraphics[width=\linewidth]{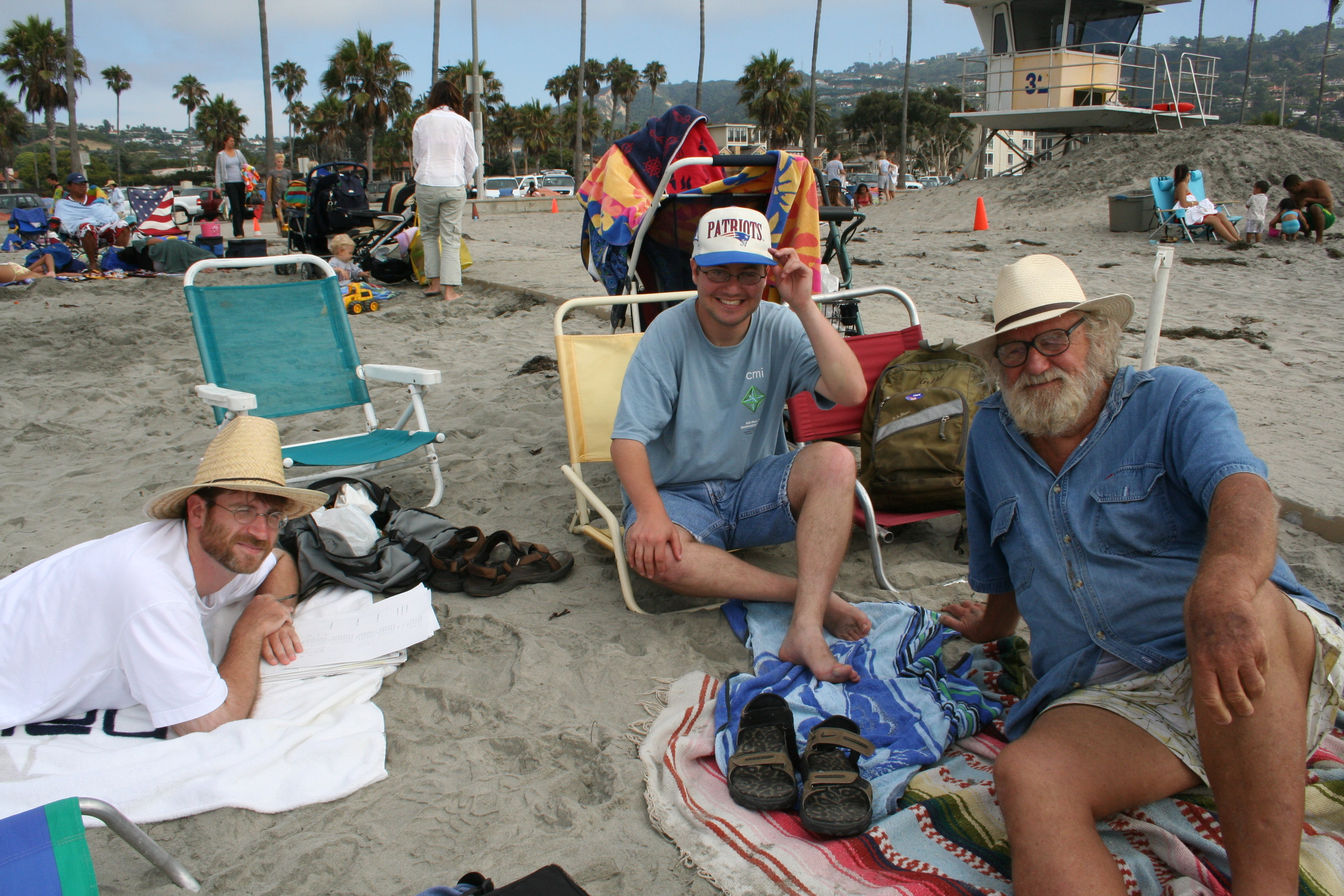}
 {\footnotesize A hallmark math session with\\
 David Little and Greg Musiker}
\end{wrapfigure}
Adriano often recalled his training as a graduate student at
Stanford doing math with Loewner in a van overlooking the ocean.
That 
influenced what became a
signature experience of almost
everyone he worked with: sharing mathematics on the beach. 
Each afternoon he sat on the shore, yellow pads of computations strewn
in the sand, drawing students and colleagues in with 
his enthusiasm for math.

He shared so much more than his passion for math.
Often, as the sun set on a math session, he invited students
and colleagues into his home for dinner.
Surrounded by the warmth of his family,
we learned to cook,
were entertained with wild tales, and discovered how to enjoy life to its fullest.

Starting in the 1960s, Adriano supervised 36
Ph.D. students over the course of six decades. 
Even after officially retiring from UCSD
in 2013, he remained research-active, advising at least five more students as an emeritus professor.
His 90th birthday was celebrated with
a thriving math conference on the beach in La Jolla.


Adriano’s work spanned multiple areas, beginning with ergodic theory and probability and ultimately finding a 
home in algebraic combinatorics. Shining throughout his career
is his insistence that deep mathematics should feel elementary 
once properly understood.

Garsia’s early work on Bernoulli convolutions uncovered deep connections
between algebraic structure and analytic regularity
in random $\beta$-expansions. For $\beta \in (1,2)$, the 
Bernoulli convolution
\[
\nu_\beta = \sum_{\epsilon_k \in \{0,1\}} 
\delta_{\sum_{k \ge 1} \epsilon_k \beta^{-k}}
\]
describes the distribution of the random series 
$\sum_{k \ge 1} \epsilon_k \beta^{-k}$. 
The partial sums
\[
S_n = \sum_{1 \le k \le n} \epsilon_k \beta^{-k}
\,,\quad\epsilon_k\in\{0,1\},
\]
typically take $2^n$ distinct values, but for certain algebraic numbers $\beta$, many coincide.

Erd\H{o}s showed that for a \emph{Pisot number},
an algebraic integer $>1$ whose conjugates all lie strictly inside the unit circle,
these overlaps are so pervasive that $\nu_\beta$ is singular. 
By contrast, Garsia studied a class of algebraic integers (now called \emph{Garsia numbers}) 
for which $\nu_\beta$ is absolutely continuous with bounded density; these have constant term $\pm 2$ 
and all conjugates strictly outside the unit circle.  

Garsia’s key insight was to measure the {collision rate} of the sums 
$S_n$ through what is now known as the \emph{Garsia entropy},
\[
H(\beta) = \lim_{n \to \infty} {H(S_n)/n},
\]
where $H(S_n)$ is the Shannon entropy of $S_n$. 
In the Pisot case, conjugates inside the unit
  \begin{wrapfigure}{l}{0.4\textwidth}
    \centering
    \includegraphics[width=\linewidth]{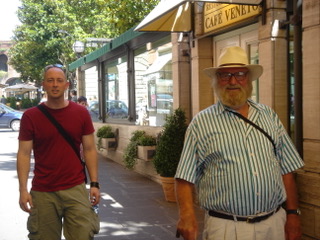}
 {\footnotesize Mike and Adriano on a Roman street}
\end{wrapfigure}
circle force systematic cancellations 
among the terms $\beta^{-k}$, producing heavy overlaps and driving the entropy 
below $\log \beta$, a signature of singularity. For Garsia numbers, the algebraic 
separation of conjugates prevents excessive overlap, keeps the entropy at its 
maximal value $\log 2$, and ensures that the measure spreads out enough to be 
absolutely continuous. 

Garsia’s entropy-based viewpoint reveals that the analytic behavior of 
these measures is governed not by chance, but by number–theoretic properties 
of $\beta$ itself. It remains central to the modern theory, especially in 
the unresolved case of Salem numbers, where the 
behavior of $\nu_\beta$ is still mysterious. 
Already in his earliest work, Adriano displayed a knack for uncovering the structural 
mechanism that makes a complicated object behave as it does.

That same instinct appears in in the {\it Garsia--Rodemich--Rumsey (GRR) inequality}. 
This inequality provides a bridge between local increment bounds and global regularity. 
Instead of relying on grueling pathwise estimates, it shows that suitable integral control 
of increments yields a quantitative modulus of continuity. 
This deterministic result underlies many regularity arguments for stochastic process sample paths.

His structural insight also guided his contributions in ergodic theory.
In that area, he is
often remembered for a three line proof of an
inequality due to Hopf, which surprised many
in the field because earlier proofs of Hopf’s result were 50 pages long.  Adriano had an ability to ``rewrite''
mathematics, stripping away clutter until the logic of a theorem 
seemed not just clear, but inevitable. 
This talent brought him true joy
and he published several books, including two coauthored with \"Omer E\u{g}ecio\u{g}lu while in his nineties.

In the mid--1970s, Adriano’s research trajectory shifted when Alain Lascoux visited UCSD and 
sold him on the beauty of combinatorics. 
At the time, combinatorics was often
looked down upon by mathematicians working
in fields of ``higher stature''
that emphasized abstraction and sweeping
generality.
Adriano felt exactly the opposite. He was drawn to combinatorics precisely because its aim is a solution requiring 
minimal background, as elementary, transparent, and human as possible.

This philosophy shaped the way he worked. After many years, he described his style by coining the portmanteau 
\emph{``manipulatorics,''} blending ``manipulation'' and ``combinatorics.'' Many of his papers prove combinatorial 
results through algebraic identities and direct calculation.

\begin{wrapfigure}{r}{0.3\textwidth}
  \centering
  \includegraphics[width=\linewidth]{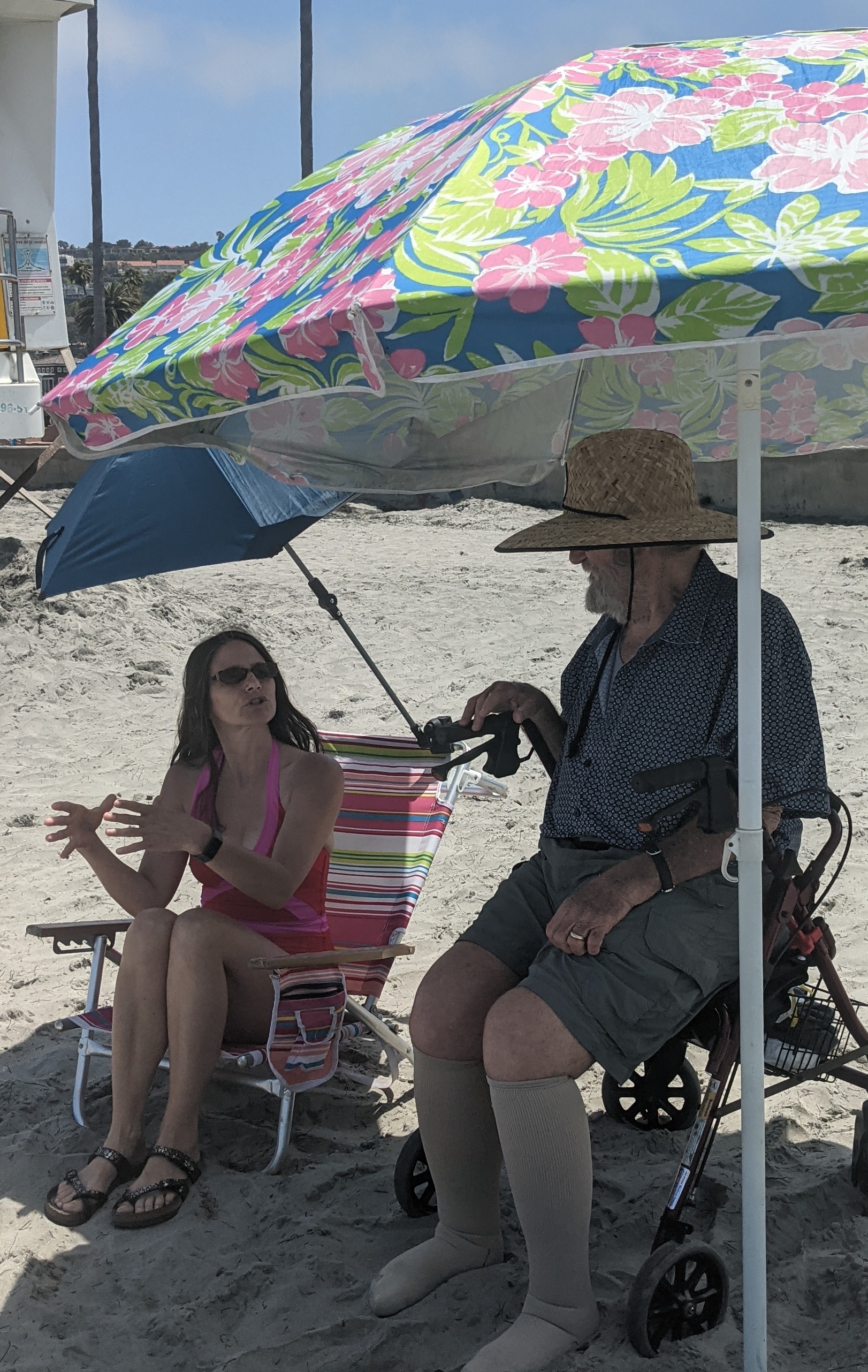}
  \caption*{\footnotesize Sharing with Jennifer his
  passion for mathematics}
\end{wrapfigure}

One striking illustration of this approach is his joint work with his student Stephen Milne on a bijective 
proof of the Rogers--Ramanujan identities, a challenge for which George Andrews had famously offered a \$50 prize. Their 
efforts led to the \emph{Garsia--Milne involution principle}, a method that reduces a signed combinatorial sum to a positive 
one by systematically combining involutions with a bijection. Rather than producing an explicit bijection outright, the 
method iterates these maps to compute a non-explicit bijection between the target sets. They collected the prize, but more 
importantly, the method itself proved more influential than the particular identity it was designed to explain.

Adriano found his intellectual home in algebraic combinatorics, where he collaborated enthusiastically on problems spanning everything from Cohen–Macaulay posets to MacMahon partition analysis. Amid this breadth of interests, there was one direction that captured his attention and remained at the forefront of his mind for the rest of his life.

Adriano invited Ian Macdonald to lecture at UCSD in 1988. There, Macdonald
explained how Kevin Kadell’s work on a generalization of Selberg’s
integral gives rise to a mysterious family of polynomials in
$\mathbb Q(q,t)[x_1,x_2,\ldots]$. Macdonald showed that these polynomials form an
orthogonal basis for the space of symmetric functions, but little else
was understood. Nevertheless, Macdonald conjectured that when expanded
in a shifted Schur basis, the coefficients $K_{\lambda\mu}(q,t)$ lie in
$\mathbb N[q,t]$. The definition was so opaque that it even obscured the
fact that the $K_{\lambda\mu}(q,t)$ are polynomials in $q$ and $t$.

Garsia was immediately intrigued and had an insight that added to
Macdonald's work, inspiring combinatorialists, algebraists, and
geometers for decades to come.  He built on the use of
$S_n$-representations via the Frobenius map,
\[
\mathcal F(V) = \sum_{\lambda}
(\text{multiplicity of } V^\lambda \text{ in } V)\, s_\lambda .
\]

Garsia proposed constructing an underlying representation-theoretic
object starting from a Schur positive symmetric function.
He introduced the {\it modified Macdonald polynomials},
\[
H_\mu(x;q,t) = \sum_{\lambda} K_{\lambda\mu}(q,t)\, s_\lambda \,,
\]
so that Macdonald’s positivity conjecture was recast as one
of ($q,t$)-Schur positivity, and he
proposed approaching
the conjecture that $K_{\lambda\mu}(q,t)\in \mathbb N[q,t]$
by realizing these polynomials as bigraded dimensions
of $S_n$-modules
$D_\mu=\bigoplus_{r,d} D_\mu^{r,d}$ such that
\begin{equation}
\label{e:FonD}
\mathcal F(D_\mu)
:= \sum_{r,d} q^r t^d\,\mathcal F(D_\mu^{r,d})
=
H_\mu(x;q,t).
\end{equation}

He started with the $q=0$ case where
the $K_{\lambda\mu}(q,t)$ reduce to the
\emph{Kostka-Foulkes polynomials}.
Garsia worked with Claudio Procesi to give a new
proof that $K_{\lambda\mu}(0,t)\in\mathbb N[t]$
by finding graded $S_n$-modules
(the {\it Garsia-Procesi modules})
whose Frobenius characteristics are $\tilde H_\mu(x;0,t)$.

Adriano then worked with Mark Haiman to
define bigraded $S_n$-modules satisfying \eqref{e:FonD}.
They defined $D_\mu$ as the subspace of the polynomial ring
${\mathbb Q}[x_1,\ldots ,x_n;y_1,\ldots,y_n]$ spanned by
all derivatives of a polynomial alternant $\Delta_\mu$.
This space is an $S_n$ module under the diagonal action
$x_i \mapsto x_{\sigma(i)}, y_i \mapsto y_{\sigma(i)}$
and conjectured that the dimension of $D_\mu$ is $n!$
and the Frobenius image is a Macdonald polynomial.
This became known as the $n!$-conjecture.


Armed with this conjecture, they expected to prove
that the $K_{\lambda\mu}(q,t)$ are in $\mathbb{N}[q,t]$
within a matter of months.
Instead, Haiman's contributions
indicated a more complex picture where
the determinants $\Delta_\mu$ are elements
of the bigraded
$S_n$ module of the {\it diagonal harmonics},
\begin{align*}
 \text{DH}_n=
 \big\{
 P(x;y)\in {\mathbb Q}[x_1,\ldots ,x_n;y_1,\ldots,y_n]
 ~~:~~
  \sum_{i=1}^n\partial_{x_i}^r\partial_{y_i}^s
 P(x;y)=0,\,\,\, \forall\,\,\,  r+s>0
 \big\}~.
\end{align*}

Fueled by computer experimentation, Garsia and Haiman amassed identities
and worked to connect the representation theory, combinatorics
and the symmetric function expressions involving Macdonald polynomials.
It was in this pursuit that Garsia and Fran\c{c}ois
Bergeron uncovered the importance of an operator
$\nabla$ on symmetric functions that takes the $\mathcal{F}(D_\mu)$
as its eigenfunctions.
They conjectured that the Frobenius image of diagonal harmonics
can be expressed compactly as
$\nabla(e_n)$ where $e_n$ is the $n$th elementary
symmetric function $\sum_{i_1<\cdots<i_n}x_{i_1}\cdots x_{i_n}$. It was later discovered that $\nabla$ is part of a broader picture connected to knot theory and
the elliptic Hall algebra.

The dimension of $\text{DH}_n$ was conjectured to be
$(n+1)^{n-1}$ while
the bigraded dimension of the
sign representation in $\text{DH}_n$ 
gave rise to a $q,t$-analog of Catalan numbers, $C_n(q,t)$, specializing to 
\[
C_n(q,1) = \sum_{\pi} q^{{\rm area}(\pi)}\,,
\]
summing  over all ${1\over n+1}{2n \choose n}$
Dyck paths $\pi$ weighted by `area,'  a familiar
statistic on $\pi$.
In 2000, Jim Haglund introduced another statistic on Dyck paths, called {\it bounce}, and conjectured that
\[
C_n(q,t) =
\sum_{\pi} q^{{\rm area}(\pi)} t^{{\rm bounce}(\pi)}\,.
\]

The $\nabla$ operator encodes all of this information
in symmetric function identities.
Garsia and Haglund worked together over the summer to give a manipulatorial proof of Haglund's conjecture, connecting symmetric function identities
with the combinatorics of Dyck paths.  Simultaneously, Mark Haiman was finalizing an algebraic geometric program that led him to establish the $n!$-conjecture and the dimension of $\text{DH}_n$.

These developments initiated a research program that
Adriano pursued for the rest of his career and
more broadly launched the still thriving area
of $q,t$-combinatorics. 
Adriano contributed his
infectious enthusiasm to the field for more than two decades.
He believed that the most profound truths were
those that could be grasped with elegance and joy.
Some of his colleagues and students
share below their memories of how he touched their
lives and shaped their careers.

\vskip .3in

\section*{H\'el\`ene Barcelo}

We have been asked to celebrate Adriano’s life, while remaining brief. How can I do so when there is so much to celebrate; by celebrating one aspect of his life at which he excelled, mentoring.

I first met Adriano in 1984 as I arrived at UCSD to begin a Ph.D. in Mathematics (logic). Diane Favreau (a Québecoise who later became his wife) had heard that a Québecoise was joining the graduate program; she tasked Adriano to find me and bring me home for dinner. Adriano and Diane not only fed me but also provided a friendship that was crucial in the first years of my life in the US and that lasted to this day. Adriano also took me under his wing. In those days mentorship was not common occurrences and without his guidance I am not sure I would have completed my study here.

Adriano’s mentorship skills were simply exceptional.  He had the extraordinary gift of detecting the potential in any students that showed a deep interest in mathematics. From then on, he would do what was needed to ensure that the student would have what they needed to succeed. If it meant spend hours one-on-one teaching intricate concepts so be it; if it meant fighting  departmental policies, let’s do it without hesitation and with fierce determination; if it meant feeding them, let’s have them over for (countless) dinners, including Thanksgiving; if it meant writing letters of recommendation let’s do it so carefully that it would take him at least a day for each letter, and he had so many students!

When I arrived in La Jolla, I certainly thought I knew everything there was to know about studying in the US. Of course I did not, and Adriano was there to help without me even having to ask!   He started by teaching me cryptography so that I could become his TA. While doing this he was also planting ideas about algebraic combinatorics. In those days, for me, combinatorics was all about arranging socks (who cares), discussing baseball (what on earth are they talking about) or still, playing poker (what is a flush!) I detested it all. There was absolutely no way I would go into combinatorics! Adriano could teach me cryptography, that was fun, but I would not become a combinatorist. He was never phased by this ``know it all'' attitude of mine. He would laugh and simply suggest that I take his algebraic combinatorics course; I did, and the rest is history as the saying goes. And it was that Adriano not only taught me beautiful mathematics, but what it meant to be a mentor: meet the students where they are, when they are ready, and without them having to ask.

Adriano was so much ahead of his time. His unconditional embracing of students in the combinatorics community is one of the reasons our field is known for being friendly and mindful of all.

I was deeply marked by my studying under Adriano, and I am as deeply thankful for his teaching, friendship, and endless generosity.

\section*{Fran\c{c}ois Bergeron}
Adriano was certainly part of my actual family, like a special uncle to my sons, as well as our mathematical one. When I say “our”, I mean our research group in Montréal. Adriano was apt at adopting places and groups with whom he enjoyed working, especially if there was also access to good food. The fact that he ended up buying a house in Montréal is significant in that regard.

\begin{center}
    \includegraphics[width=.45\linewidth]{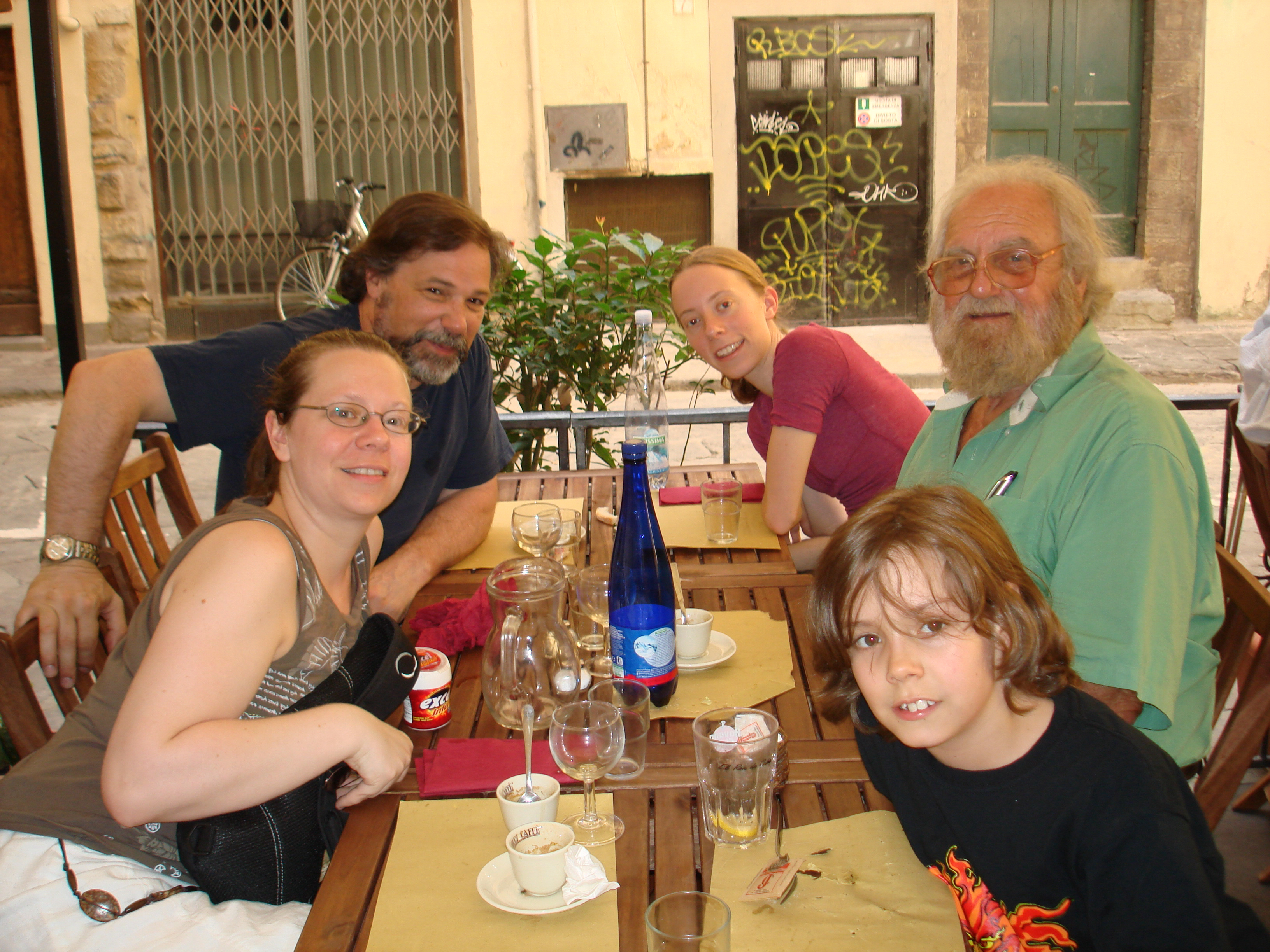}
    \includegraphics[width=.45\linewidth]{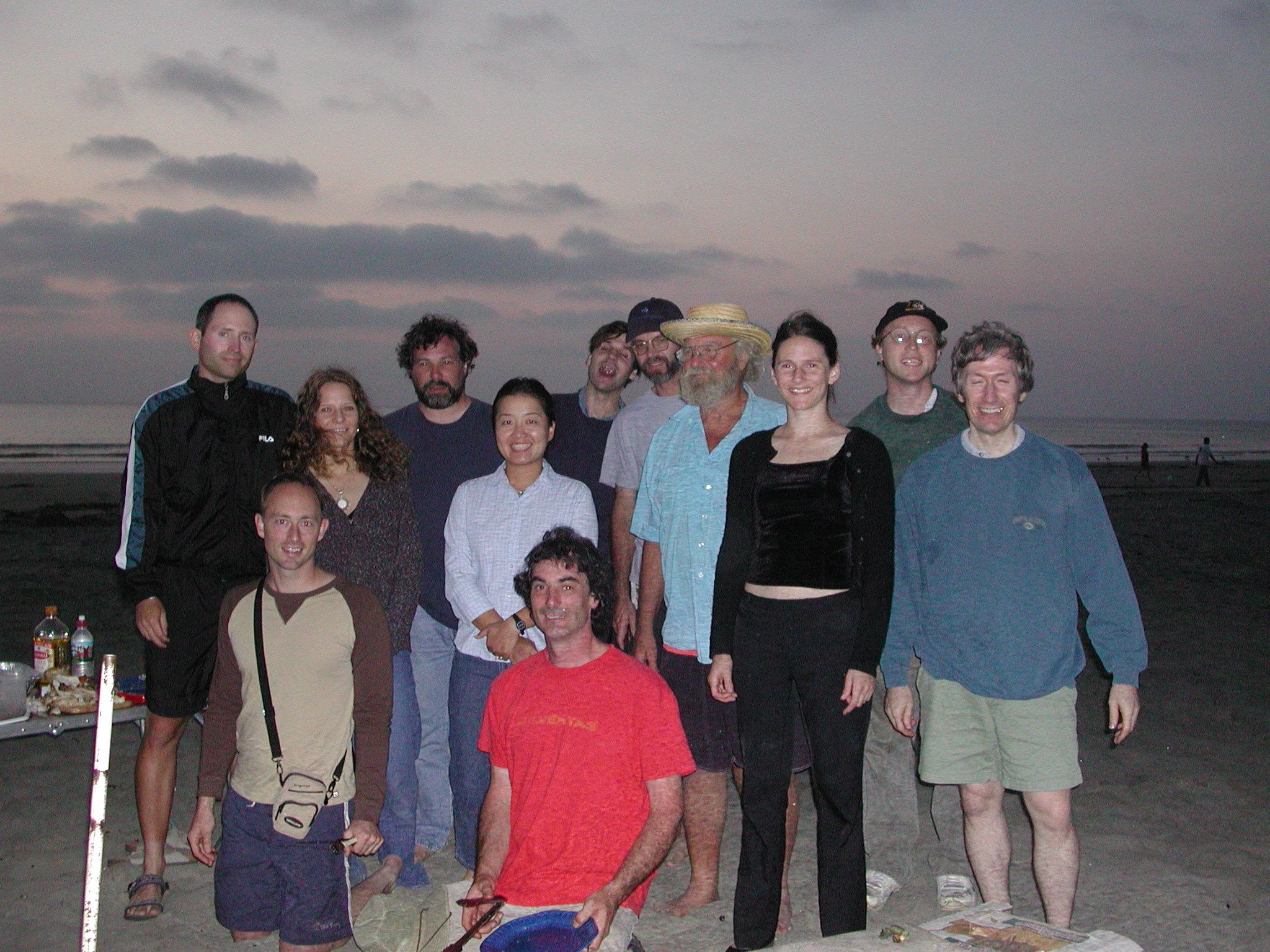}
\end{center}

I first met Adriano around 1980, when he came to visit and gave us a talk about the use of generating series in combinatorics. He was invited by Pierre Leroux who was just back from a sabbatical in San Diego, having been exposed to impressive combinatorics by Adriano and his visitors. One open question that he mentioned in this talk was to give a combinatorial interpretation for the substitution of exponential generating functions.
 This poked the interest of André Joyal, who came up with an answer in just a few days. What was striking for many of us, was the use of Category Theory in this, with combinatorial construction described by functors. At the time, Montréal was a sort of Mecca of Category Theory, and that was our usual playground. On my part, I was working on my doctoral project, having to do with a Topos Theory and abstract constructions related to Differential Geometry. I essentially knew nothing about combinatorics. 

Our group, with Pierre Leroux, Gilbert Labelle, Jacques Labelle and Hélène Décoste, started playing with Joyal’s notion of Species of Structures, learning combinatorics on the way. Thus Adriano’s talk turned out to be key in the birth of research center, going on strong almost 50 years latex. 

\begin{wrapfigure}{l}{0.5\textwidth}
    \includegraphics[width=.95\linewidth]{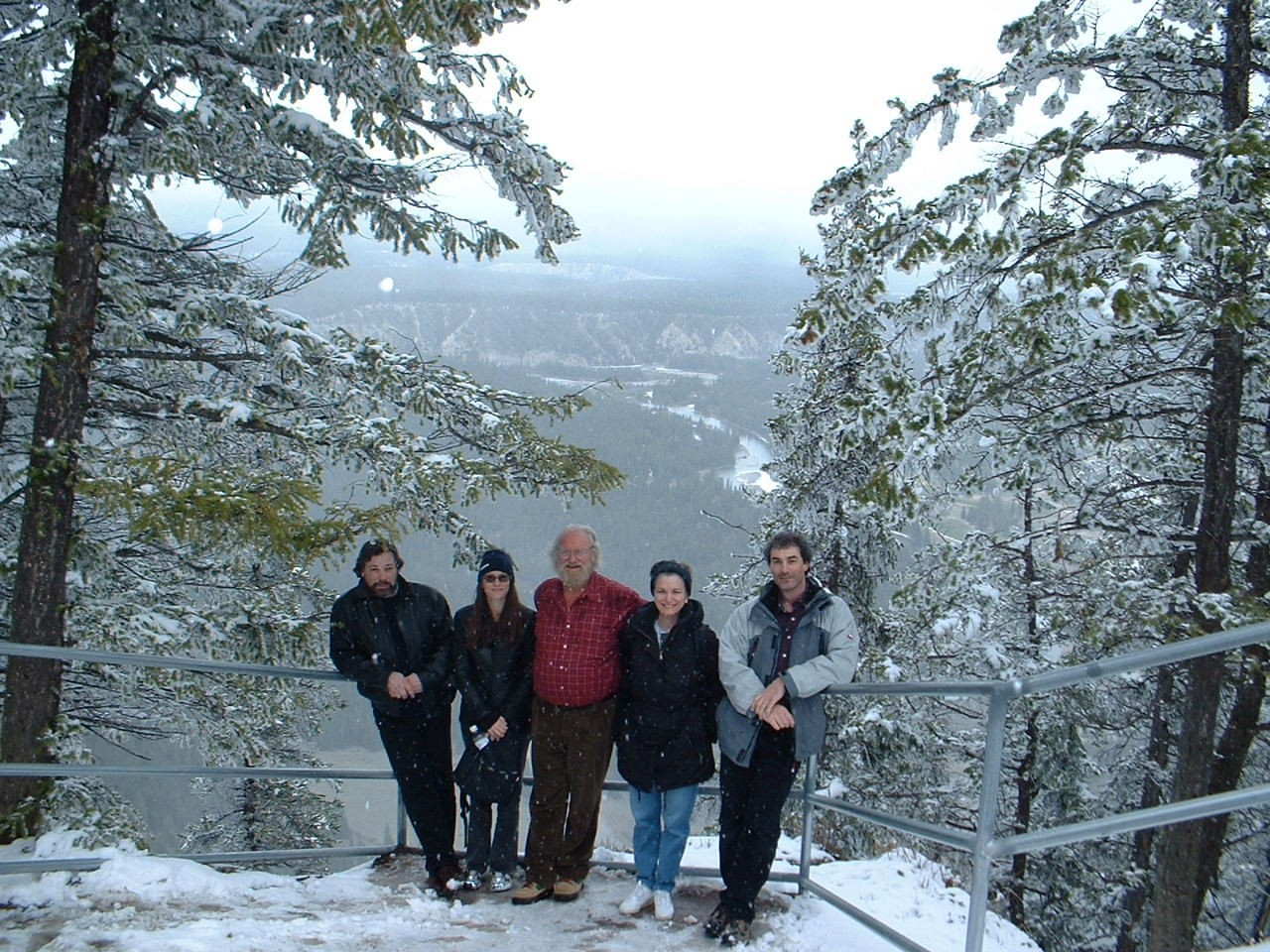}
\end{wrapfigure}

The interplay between Combinatorics and Algebra so fascinated me that I changed the subject of my thesis. From time to time Adriano would visit, or participate in some large meetings that we organized in Montréal. His beautiful combinatorics so fascinated my brother Nantel, that he went on to work with him in San Diego. You can learn about this in his contribution. What is fun about this, is that Adriano visited again a bit after Nantel’s move to San Diego.
This is the first time that we collaborated. He went back to San Diego in the middle of this collaboration, and kept on working there with Nantel. This ended up with my first joint paper, together with Adriano and Nantel, without me talking directly to Nantel. So Adriano was the bridge between us, at a time when the web was not actually in place (before 1990).

Over the years, it became an habit for me to visit San Diego in the Spring, together with my family, and enjoy Adriano’s passion for mathematics and great food. In turn, Adriano would often visit Montréal, bringing students and family.
It was always exciting to propose new ideas to Adriano, because of his amazingly enthusiastic reactions. This is certainly one of the reasons that he is among those that I have the largest number of collaborations, not only because of the beautiful mathematics, but also because of the great fun we had while doing it.

\section*{Nantel Bergeron}

\begin{wrapfigure}{r}{0.35\textwidth}
    \includegraphics[width=.95\linewidth]{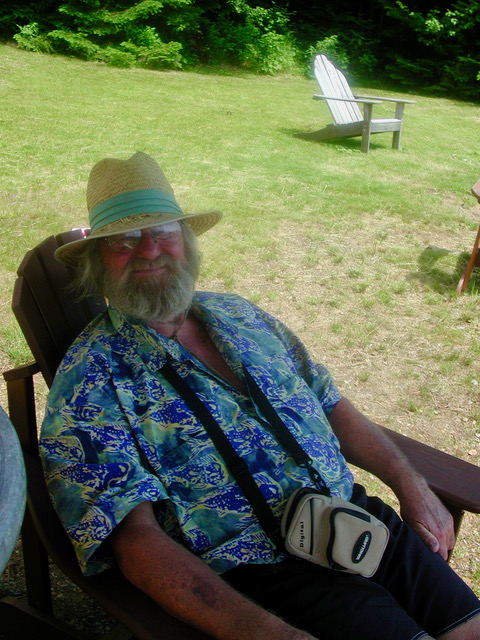}
{\footnotesize Adriano at at the lake near Mont Tremblant, Qu\'ebec.  }
\end{wrapfigure}

Of course, when thinking about Adriano, countless dinner stories come to mind. He was passionate about life, and it was around the making and enjoyment of food that one got to know him more deeply. From stories of his childhood to tales about mathematics, Adriano was an open book filled with fascinating plots and anecdotes.

I always felt that Adriano was part of my family--and that I was part of his. 

Slowly, through these wonderful dinners, I was introduced to some of the most influential mathematicians of the time! I still remember the day Sch\"utzenberger visited UCSD and Adriano wanted to cook dinner for the master. At the time, he was in transition between houses, so he used my apartment to cook. I was incredibly intimidated to host Sch\"utzenberger—this was quite something for me. And I must say, it was the first time I felt that Adriano himself was also a little bit intimidated. Perhaps he even expressed it more than he felt, just to make me feel at ease. I wouldn’t be surprised—this would be so in line with his nature.

Yes, this was a wonderful life! But somehow, mathematics was always there—woven into everything. It is through mathematics that I discovered a very different mind. I remember him going through pages of computations. Sometimes I would follow line by line, only realizing at the end the scope of what we had done. I would always marvel at that moment--``How did he do it?'' And again, he would make me feel comfortable, point to something I had done that he found wonderful.

I saw him get so excited about what those around him had accomplished, especially us—his students. In fact, that’s how he recruited me to study in San Diego in 1987. We had worked on some small representation theory problems in Montreal, and since I had solved a minor question, he got all excited and convinced me to give up my plans to go into physics and instead study with him at UCSD.

I am sure many of his students will recognize themselves in these stories.

Thank You Adriano

\newpage

\section*{Sara Billey}

\noindent
{\bf Remembrances and Lessons from Adriano Garsia}

\begin{center}
\includegraphics[width=3in]{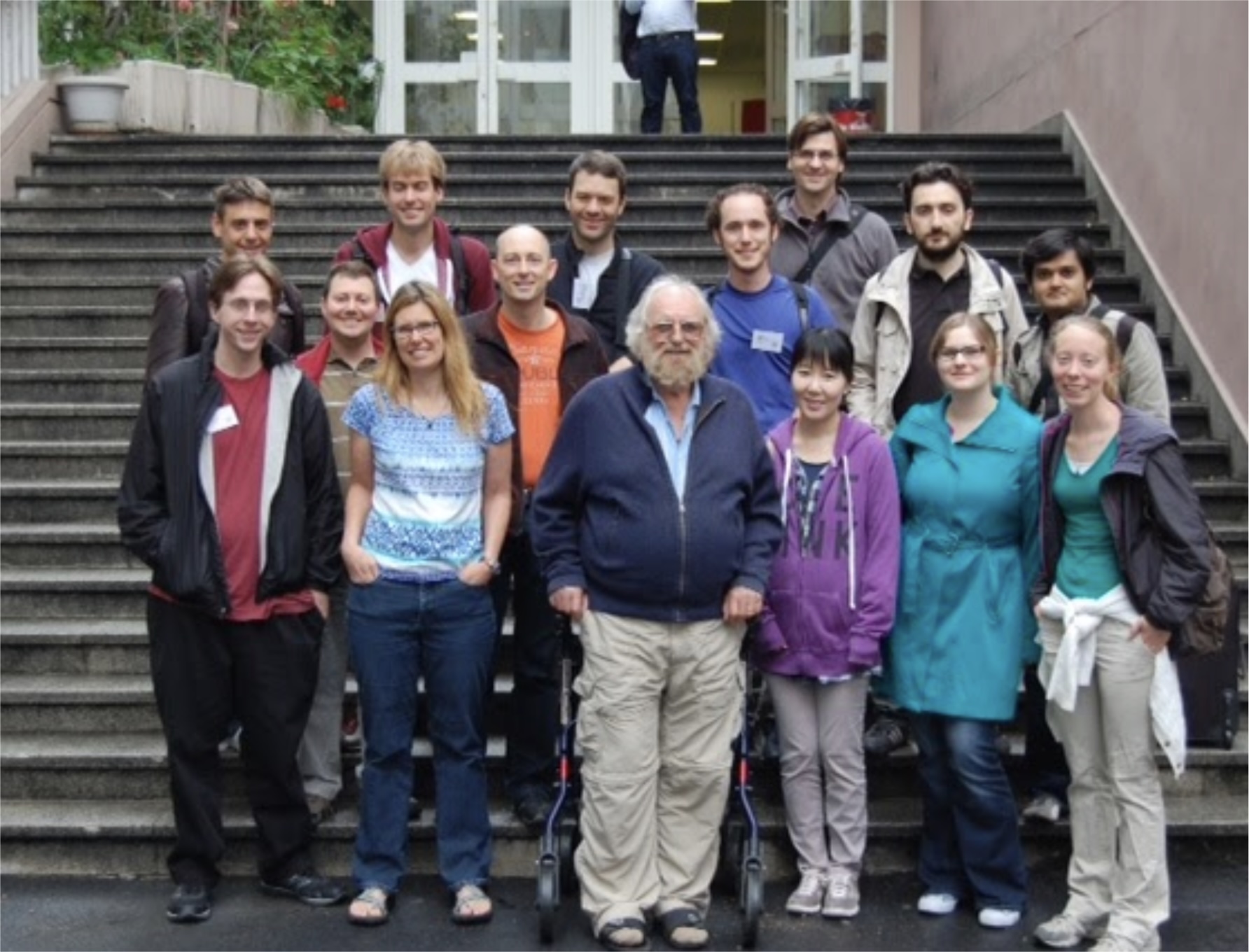}\\
{\footnotesize  Academic family photo with 3 generations!}
\end{center}

{\bf Be enthusiastic}. That was advice I heard recently in a talk for young
professionals. It was excellent advice, but it was advice I had
already learned from Adriano Garsia. Adriano brought intense
enthusiasm to mathematics, and the field was better for it. He brought
enthusiasm to food, and we cooked more passionately. He brought
enthusiasm to the beach, and we went to the beach. He brought
enthusiasm to building community, and his community grew.

Adriano valued the people around him: his family, his friends, his
collaborators worldwide, and his students. He lifted people up through
praise, attention, and genuine appreciation for what each person
contributed.

\vskip .1in
{\bf Novel contributions}. This past spring I had the pleasure of teaching
Modern Symmetric Functions. We focused on work from the past thirty
years, and Adriano’s influence was everywhere: the Garsia–Haiman
modules, the $n!$ conjecture and theorem, and the diagonal coinvariant
ring and the space of diagonal harmonics. What an extraordinary
mathematical journey he helped set in motion, beginning with his work
with Claudio Procesi and Mark Haiman and continuing through work with
Mike Zabrocki, Jennifer Morse, Jim Haglund, Glenn Tesler, and many
others. Adriano inspired generations of mathematicians to work on
these challenging problems and to bring to them tools from algebraic
combinatorics, representation theory, symmetric functions, algebraic
geometry, and beyond.

In a sense, the story of the $n!$ Conjecture began even before the
1990s, in Adriano’s work on the coinvariant algebra and related
quotients. The coinvariant algebra of $S_n$ is the quotient of the
polynomial ring in n variables by the ideal generated by the
positive-degree symmetric polynomials. It is a graded $S_n$-module under
the action that permutes the variables and, after forgetting the
grading, is isomorphic to the regular representation of $S_n$. A natural
question remains: which basis best reflects this action?

The descent monomial basis, now often called the Garsia--Stanton
descent basis, became a central candidate. In the early 1980s, Garsia
and Stanton carefully studied bases for related quotient rings and
analogues of the Vandermonde determinant. This in turn inspired Ed
Allen’s work on straightening algorithms for descent
monomials and related quotients. Later work connected the
Garsia–Stanton basis to descent representations and to the
Lusztig--Stanley description of the graded decomposition of the
coinvariant algebra via the major index statistic on standard Young
tableaux. These threads helped prepare the ground for the diagonal
world that made the n! conjecture so compelling.

\vskip .1in
{\bf Build knowledge by bringing people together}. Adriano loved having
mathematical visitors come to San Diego to share their research with
his group. During my years at UCSD, we had visits from Richard
Stanley, Bruce Sagan, Claudio Procesi, Michelle Wachs, and many
others.  Ian Macdonald was a regular visitor in the 1980s and 1990s,
and he often gave courses during his winter visits. Two of those
courses had an enormous impact on me, and I remain grateful to Adriano
for creating the opportunity.

First, Macdonald invested a great deal of time in understanding the
work of Lascoux and Schützenberger on Schubert polynomials. In his
course, he made the theory accessible to us, supplying a coherent and
rigorous development from the original papers. The resulting Notes on
Schubert Polynomials were typed up by UCSD graduate students at the
time, including Nantel Bergeron, Hélène Barcelo, Arun Ram, and Tamsen
(Whitehead) McGinley. Those notes became a foundation for many of us
learning combinatorial Schubert calculus, and they inspired me and
many others in the field. They are now available online from LaCIM.

Macdonald also gave a second course on root systems during my first year of
graduate school. His handwritten lecture notes were as complete and
polished as a published book, with theorem numbers, full proofs, and
citations. 

\vskip .1in
{\bf The culture matters}.  The research environment Adriano created was essential to the
success of his students and to the strength of the wider algebraic
combinatorics community.  His enthusiasm continues to inspire mathematicians 
world wide.   And, his favorite beach, La Jolla Shores, is a lot lonelier without him.

\section*{Anders Bj\"orner}

Adriano Garsia was my friend. I was proud that he often referred to me as his ”mathematical brother”. 
 
Adriano started research in analysis in the 50s, as a student of Charles Loewner and Marcel Riesz, who were both at Stanford at that time. His proof of the mean ergodic theorem from 1955 is considered by many the best way to understand that theorem. In the 70s, Adriano moved on to combinatorics, more precisely algebraic combinatorics, where he soon was recognized as one of the leaders in the field.
 
Adriano and I started talking about combinatorics when we first met at an AMS general math conference in Columbus, Ohio, in 1978. His enthusiasm for combinatorics was contagious, and our conversation continued at UCSD a few months later – what was to be one of many visits. His favorite ”seminar room” was under a parasol at the beach in Del Mar, where even the most intellectual Parisian mathematician succumbed to balancing his notebook on a towel filled with sand. I especially appreciated that coffee breaks were replaced by a pause to ride the waves on a boogieboard! This is an excellent way to clear one’s mind – I can recommend trying it.
 
Adriano’s beach seminars produced beautiful results. Among them is his work in 1977 (with Michelle Wachs) on optimal binary trees, and the proof in 1993 (with Mark Haiman) of the n! conjecture. Of lesser importance but dear to my heart, was the expository paper on Cohen-Macaulay posets that I wrote with Adriano and Richard Stanley in 1982. 
 
While we can acknowledge Adriano’s contributions to mathematics, how does one capture in a few words the enthusiastic spirit of this man?
 
Adriano’s students bear witness to the fact that he was a very inspiring teacher, always full of ideas, always asking questions.  His enthusiastic style of mentoring encouraged everyone to believe in their abilities. The gratitude of Adriano’s students will be his legacy.
 
Adriano was a first-class story teller. I remember a colloquium dinner at UCSD where Saul Bellow was present. People were in awe to be in the company of a Nobel laureate. But also Bellow himself had reason to be impressed.  The next day he thanked Adriano profusely for the captivating stories he had shared about his childhood in Tunisia and adolescence in Rome during the second world war.
 
Adriano showered his family and friends with his generosity and kindness. On one visit to UCSD, he handcrafted the wood furniture for an unfurnished house he’d rented for my family in Del Mar. And he helped us evacuate this house and move in with him when there was a tsunami warning! 
 
Finally, let me say that Adriano was at the core a passionate and proud Italian. This flavored his life and his cooking (his secret was adding handfuls of salt to the spaghetti water!).  Only espressos and speedos, never American coffee or beachwear.
 
I will remember Adriano not only as a mathematician, but for who he was as a person. Adriano was a huge figure in my life. I miss him dearly.

\section*{Tullio Ceccherini-Silberstein}

\begin{wrapfigure}{r}{0.5\textwidth}
    \includegraphics[width=.95\linewidth]{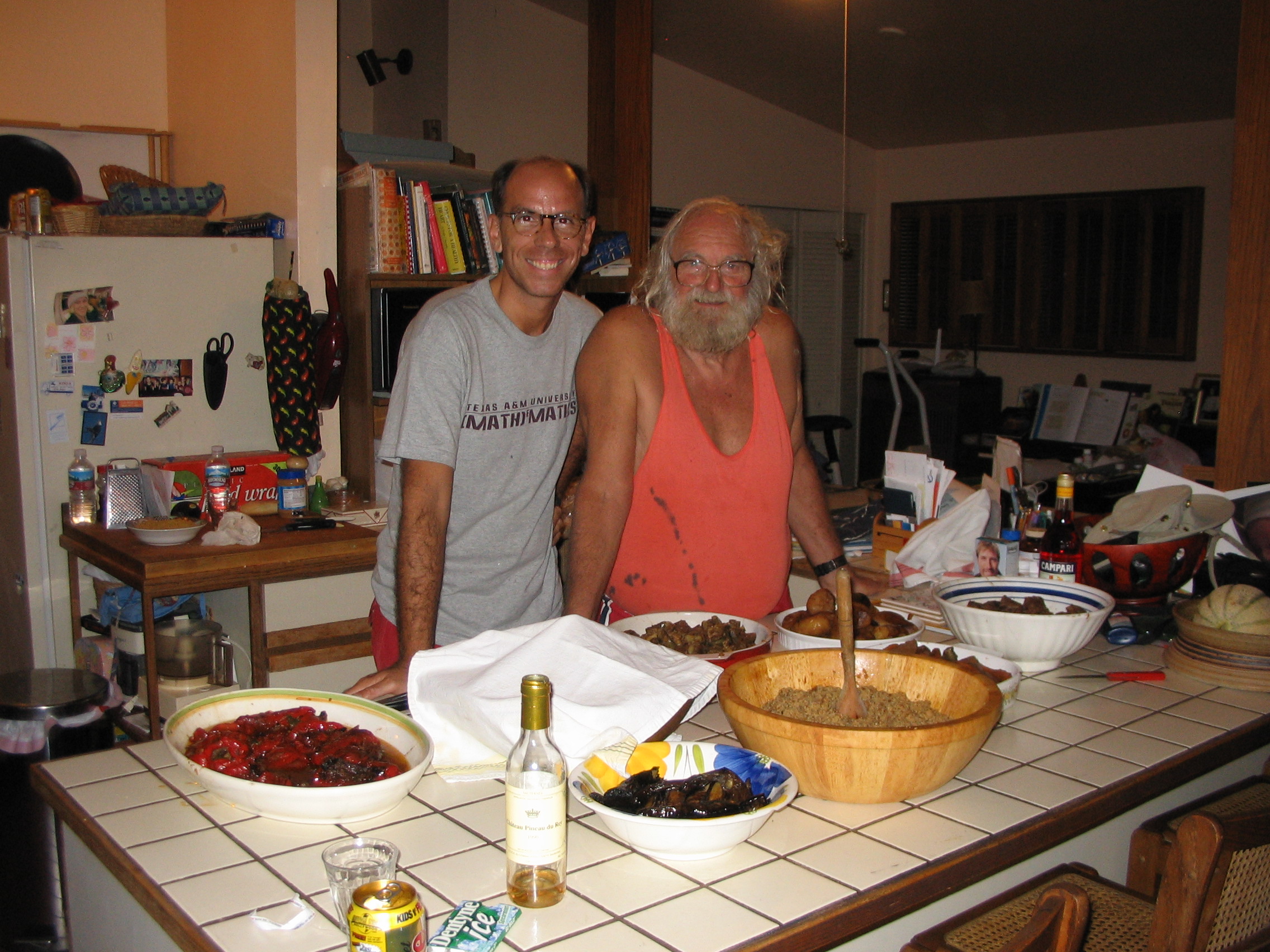}
{\footnotesize With Adriano in his kitchen, after preparing a couscous, circa 2006}
\end{wrapfigure}

It is both a pleasure and a great honor to contribute a few personal recollections about Adriano.
Given Adriano's extraordinary personality, his remarkably broad interests, and the uncountable
adventures that marked his long and fascinating life, one could easily fill several volumes with
anecdotes about him. I shall therefore resist the temptation to recount his youth in his native country,
Tunisia (which, as he liked to joke, made him an exceptional ``African-American'') -- where his
grandfather used to serenade his lover, the wife of a Russian admiral stationed in Tunis, by singing
arias from Mussorgsky's Boris Godunov -- or his early romantic adventures, his family's move to Italy
after the war (including the story of his beautiful cousin, who turned down the marriage proposal of a
charming and wealthy young man of Italian descent, the owner of a prosperous agricultural estate,
because she could not bear the thought of living on a farm, far from the bustling city center of Tunis,
with its elegant shops and fine restaurants -- only to spend the rest of her life, after the Second World
War, in Velletri, a tiny town some 40 km south of Rome), his passion for cooking (especially his
unforgettable couscous, above all the fish version known as ``agliotta''), and many other stories that
made every conversation with him both delightful and unforgettable.

I will only recollect three episodes of a mathematical interests.

I first met Adriano in Rome in January 2003, when, at the invitation of Claudio Procesi, he gave a
lecture at the University of Rome ``La Sapienza.''
Although the topic was Algebraic Combinatorics -- a
subject about which I then knew virtually nothing and in which I had not yet developed any interest --
I was curious to hear the mathematician whose proof of the Ergodic Theorem had become famous for
its remarkable brevity and elegance. A week later, I travelled to UCSD to visit Efim Zelmanov for a few
days, and discovered that Adriano -- who had meanwhile returned to San Diego -- was a professor
there. I wrote to him, and he immediately invited me to dinner. We agreed to meet late in the
afternoon at La Jolla Shores, on the beach overlooking the Pacific Ocean, where he was working with
Richard Stanley. When I arrived, I immediately recognized the two mathematicians, deep in discussion
under a beach umbrella. Introducing myself to Richard, I could not resist saying:
``Dr. Stanley, I
presume?''
Later that evening, the three of us enjoyed a wonderful dinner at Adriano's home on Mt.
Armet Drive, where I also met his wife Diane and their daughter Gabriella. That evening marked the
beginning of a beautiful friendship.

In the late 1940s, as a mathematics student at the University of Rome ``La Sapienza,'' Adriano
attended the Mathematical Analysis courses taught by Mauro Picone and his young assistant Gaetano
Fichera. Having recognized his exceptional talent, they played a decisive role in enabling him to
pursue graduate studies in the United States.

In 1988, Fichera -- then in his late sixties -- taught me a beautiful course in Higher Analysis devoted to
Fourier series. One of its highlights was Marcel Riesz's celebrated theorem on conjugate series.
Fichera had the delightful habit of introducing mathematical topics with personal recollections and
historical anecdotes. On that occasion, he told us that he had known Marcel Riesz personally and had
enjoyed his friendship. Many years later, Adriano told me that when he was a graduate student at
Stanford, his advisor Charles Loewner asked him to look after Marcel Riesz during a several-month
stay at the university. Every day Adriano would pick him up at his hotel, drive him to the university,
and take him back in the evening. Marcel Riesz was not only a great mathematician but also a
remarkably charming and generous man. He frequently invited Adriano to breakfast, lunch, or dinner,
and their conversations ranged far beyond mathematics, although mathematics naturally occupied a
privileged place. I found it deeply moving that one of Adriano's earliest mathematical passions was
precisely the theory of trigonometric series. He later wrote the beautiful booklet ``Topics in Almost
Everywhere Convergence,'' containing, among other things, the elegant proof of the Ergodic Theorem
that I mentioned above. Years later, he presented me with a copy, inscribed with a warm personal
dedication -- a gift that I have always treasured.

In the summer of 2011, I taught in the Summer Session at UCSD and brought my two sons,
Giacomo (who has since become a mathematician!) and Tommaso -- then thirteen and ten years old --
with me to San Diego. One evening, while we were all enjoying a delicious dinner prepared by
Adriano
(with our enthusiastic, if not particularly useful, assistance), he treated us to a truly Mark Twain-like
episode from his own life: he read aloud his newly published obituary from the online Notices of the
American Mathematical Society! The story was as bizarre as it was amusing. Some time earlier, a
woman named Adriana Garcia had been murdered by her husband in Santa Ana. Someone, hearing
the news on the radio, misunderstood the name, concluded that it was Adriano Garsia who had died,
and informed the AMS. The result was that Adriano had the rare privilege of reading his own obituary
while still very much alive.

By then, Adriano had already told us many stories about the many remarkable mathematicians he had
known personally. Shortly after completing his Ph.D., he spent a year working for the U.S. Army on a
classified project concerning the optimal shape of a shell containing a nuclear device -- a problem
involving Riemannian geometry and potential analysis. Although the work remained confidential, it
brought him into contact with John Nash, who had recently proved the celebrated embedding
theorems that now bear his name. Around the same time, through their common interest in
trigonometric series, Adriano also became friends with Paul Cohen, several years before Cohen's
celebrated independence proof for the Continuum Hypothesis. Sadly, when Nash's mental illness
began to manifest itself, his physicians discouraged friends and colleagues, including Adriano and
Cohen, from maintaining contact with him.

Recall the famous scene in the movie ``A Beautiful Mind'' in which Nash finally realizes that he is
suffering from schizophrenia. He understands that the little girl with whom he has been speaking for
years cannot be real because, as time passes, she never grows up. At one point she tugs at his sleeve,
and Nash abruptly pushes her away, exclaiming,
``Go away! You don't exist!''
Returning to our delightful dinner, my sons and I naturally asked Adriano whether he had ever
reestablished contact with Nash.
``No,'' he replied. I could not resist imagining the following scene.
Adriano meets Nash after many years and says,
``John, it's Adriano Garsia. Do you remember me?''
Nash looks at him in disbelief and answers,
``You can't be Adriano. I read his obituary in the Notices of
the AMS. You don't exist!'' Adriano laughed heartily, and so did the rest of us.

Those who had the privilege of knowing Adriano will remember not only a brilliant mathematician,
but also a man of extraordinary generosity, curiosity, a wonderful sense of humour, and deep
humanity. Every meeting with him was an adventure, every conversation a lesson, and every dinner a
celebration of friendship and life. Thank you, Adriano. I miss you very much.

\section{Michele D'Adderio}

When I arrived in UCSD to do my PhD in algebra, I knew Adriano Garsia by fame. So, after a few days of settling in, I sent him an email:

\smallskip 

\verb|”Dear professor Garsia,|\\
\indent \verb|this is Michele D’Adderio, a new PhD student at UCSD.|\\
\verb|I did my master thesis with professor Procesi,[...] I wanted to have|\\
\verb|more information about your course ’Topics in combinatorics’...”|

\smallskip 

Next day I got an email from Garsia:

\smallskip 

\verb|”Caro Michele, incontriamoci domani mattina a La Jolla shores.|\\
\verb|Io saro' vicino alla tower 32. Questo e' il mio numero di telefono: 619...”|

\smallskip

\begin{wrapfigure}{l}{0.5\textwidth}
	\includegraphics[width=0.95\linewidth,clip=true,trim=25mm 15mm 20mm 25mm]{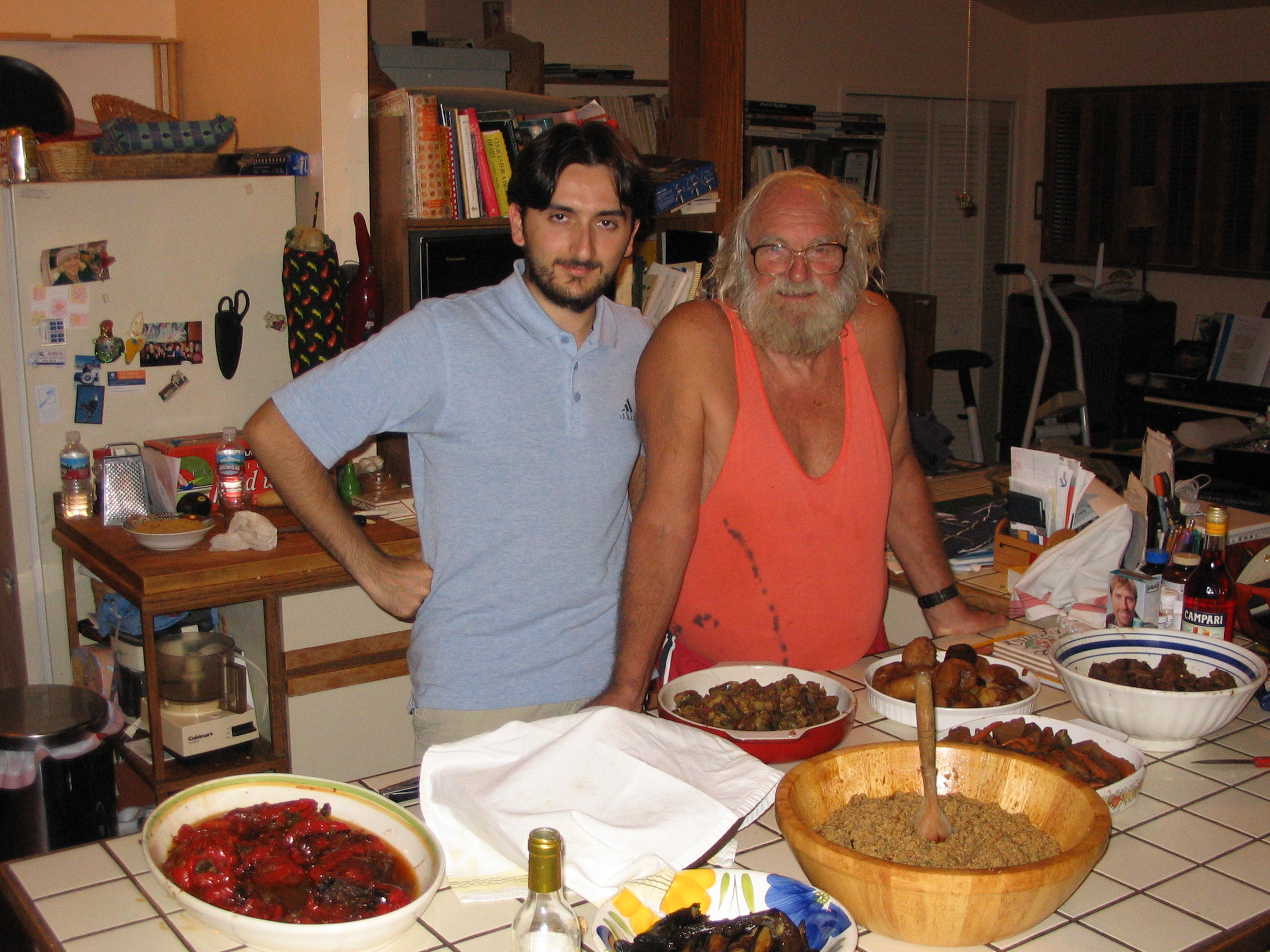}
	\caption{First dinner with Adriano.}
	\label{Fig:m_ary_tree}
\end{wrapfigure}

(``Dear Michele, let us meet tomorrow morning at La Jolla shores. I will be next to the lifeguard tower 32. Here is my cell phone number: 619...'')

\smallskip

The next day I went there, and that was when I first met Adriano, right in his natural habitat: the beach! And that was also when I got my first lecture on the representation theory of the symmetric group, symmetric functions, Frobenius characteristic, etc.

\smallskip

A few days later, the beginning-of-the-year drink of the math department took place.
There I saw Adriano talking (in Italian) with Tullio Ceccherini-Silberstein (whom I met in that occasion): he called me, he introduced me to Tullio, and he invited both of us to have dinner at his place the next day.

So there I was, welcomed in his home, helping him preparing his famous couscous, having dinner together with his wife Diane, his daughter Gabriella, and Tullio. And all this just after we talked literally twice!

\smallskip

In the four years I spent in San Diego I met Adriano regularly, several times per week. 
He treated me fully as one of his students, though I have not been formally one of them (I did my PhD in noncommutative ring theory). The life changing experience of being a PhD student of Adriano can be partially grasped by some of the stories and thoughts found in this article. But it is hard to convey in its fullness.

During those years and the ones afterwards I witnessed and heard from him innumerable anecdotes and adventures, the one I reported above being not atypical. E.g.\ once I went to Kauai with Adriano and his PhD student at the time, Angela Hicks, to do mathematics (you can read Angela's account in this article). Another time I accompanied him in his trip to France to receive his honorary doctorate from the University of Bordeaux (he found ironic to sign the ``book'' of recipients, whose first signature was from Charles De Gaulle, his ``enemy'': Adriano's family farm in Tunisia was expropriated by the French after the second world war, forcing them to flee to Rome). Needless to say, those trips where memorable.

\smallskip

Everyone meeting Adriano for the first time, at a dinner, would be forever impressed by his entrancing charisma, his playful wit, his youthful enthusiasm, and his operatic storytelling.
But when it came to mathematics, his allure would turn into an irresistible whirlpool of excitement.

His influence on my mathematics has been immense. I spent countless hours learning mathematics from him: every time we met he had something ``fantastic'' to show me. Like an apprentice in the workshop of a Renaissance master, I watched him working with his yellow notepad on the beach, or programming at the computer in his office. I have always been both fascinated and intimidated by his unique way of looking at mathematics: he would reduce every problem to its most elementary essence, and then attack it with almost bare hands. His mastery of symmetric functions identities (an example of what he called ``manipulatorics'') was the skill that scared me the most. 

After I finished my PhD in algebra, I moved to combinatorics, but not the one Adriano has been teaching me all those years. It took me a few more years to pick up one of the open problems that he really cared about, and to finally start working on it. Ironically, I ended up (re)learning all his scary manipulatorics of symmetric function identities from his articles, and I am still joyfully struggling with it to these days. 

When I communicated to Adriano that I started to work on the problems he liked, he quoted an old joke to make a statement of the effect ``You are finally home!''. 

In 2019, at Garsiafest, the conference dedicated to his 91st birthday, he would loudly say that he finally converted me to his favorite mathematics. He (half) joked that my ``conversion'' was his ``revenge'' against the algebraists who ``stole'' Jason Bell from algebraic combinatorics (Jason started his PhD with Adriano, but he ended up doing his thesis with Lance Small).

At that same conference, a special dinner with Adriano and all his PhD students took place. The day after Adriano told me about the dinner, how the evening had gone, and how everything and everyone was ``fantastic''. In that conversation he invited me to have dinner at his place, with his family, that same night. I believe Adriano invited me to compensate the fact that his ``adopted student" was not there the night before: once again he made me feel part of his family in every sense of the word. That fabulous dinner was also the last time I talked to Adriano in person.

In the following years I sporadically talked to him on the phone. It was lovely to hear his voice, his stories. In our conversations he has always done most of the talking; but in his last years the percentage of his talking time was getting even closer to 100\%. Of course he was fully aware of this phenomenon, and after a long stretch of him continuously talking, he would say with his usual wit ``Sorry Michele, these days I am in \textit{download mode}...''.

I miss him a lot, and I think about him very often. When I see a nice theorem or a beautiful proof I often think ``Adriano would have loved this one...''. I talk about him a lot to my own students, and it always feels awkward that I have to refer to him as ``Garsia''... while in my mind he is simply ``Adriano''. 

\section*{Mark Haiman}
I got to know Adriano when I was seeking my first permanent position after being a postdoc.  He ended up recruiting me to UCSD, where I stayed for 10 years before moving to Berkeley.  That period, and after, was for me a time of wonderfully fruitful collaboration with Adriano and others in his circle.

Adriano and I both attended the combinatorics program at the Mittag-Leffler institute during the winter quarter of 1992, my first year at UCSD.   Adriano's wife Diane and young daughter Gabriella were there with him.  After a while, when the grey Stockholm weather (so far from what he grew up with in Tunisia!) had begun to drive Adriano crazy, he announced that he was taking his family to Rome for a week, and that I should come too.  That was when he introduced me to Claudio Procesi, and taught me to drink espresso.  It was the first of several trips to Italy with Adriano and sometimes with Diane and Gabriella as well.  

Another time, after a conference in honor of Adriano's 70th birthday in Taormina (at a hotel on the beach, of course), the four of us made a driving tour the whole length of Italy and on to Avignon and the TGV to our flights home.  Adriano loved company, knew and loved Italy, and loved showing places he loved to others: Sicily and Taormina, Sorrento and the Amalfi coast, Pompeii, Rome, Pisa, Fiesole, the Cinque Terre (this last was new to him, as I recall, and all of it was new to me).  The whole excursion was a fabulous experience.  We got on well together, with just one quirk---both Adriano and Diane were prone to car sickness.  Things were fine when Adriano drove, as the driver doesn't feel it, and Adriano had learned to make adjustments to keep Diane happy, but naturally I needed to share the driving, especially since Diane was uncomfortable driving on Italian roads.  So I was given two instructions: (1) take all curves very slowly, and (2) never touch the brakes!

It would be impossible to relate here all of Adriano's extraordinary talents.  Of course, everyone who knew him remembers the long conversations over amazing meals that he would cook.  He was also a maker of other things physical as well as mathematical.  His study upstais at the house in La Jolla was filled with bookcases and a wrap-around desk of his own construction; later he built a large gazebo in the backyard for eating outdoors, hiring day laborers to help with the work, and communicating with them in some hybrid of their Spanish and his Italian.  The seminar room in the UCSD math department having had insufficient blackboards originally, he built two new large ones on sturdy movable stands, which were still in use the last I knew.

As personalities and collaborators, he and I were an odd match: Adriano was the epitome of gregariousness, while I tend to be solitary.  Whenever the weather was good, he would set himself up in the afternoon at La Jolla shores or Del Mar beach, with his towels, folding chair, umbrella, bag full of math, and students; I didn't like to work in the sand, and would only join him there occasionally.  The way Adriano would describe our working style, he was the optimist and I was the pessimist.  He would come up with one idea after another, thrilled every time to believe that this new idea must be the key to solving whatever we were trying to solve; I would then search for an example in which it wouldn't work.  This however turned out to be an excellent way to collaborate, for when you open many doors, then shut the ones that lead to blind alleys, you eventually find the door that leads to the place you want to get to.

\section*{Angela Hicks}

When I became Adriano's graduate student, there had been a significant gap since his last student, and rumors were swirling that he was no longer advising. Ha! There were at least four after me, and he continued mentoring for nearly a decade and a half. Before I signed on, a senior student of Jeff Remmel warned me that she had no interest in ``drinking water from a firehose.'' I quickly understood the analogy, but I benefited immensely from the absolute deluge of mathematical attention Adriano showered on us.

As an advisor today, I try to emulate his boundless enthusiasm and open-door policy, although I meet with my students slightly less than the one-or-more times a day Adriano demanded mathematical companionship. With some grads, he was highly competitive, racing them to correct answers. I ducked that by taking on the purely combinatorial parts of our arguments, leaving him to dominate the equations and ``manipulatorics'' that were his passion. It made us great collaborators, though I had to spend a lot of time consulting his papers once I began working independently.

Adriano was a beautiful writer who prided himself on checking every single line of my arguments (including the few he slept through and I had to repeat). However, he had far less patience for student writing. When I was preparing my thesis, he became incredibly grouchy that writing was stealing me away from our daily sessions. When he ``finally'' received my first draft, he read it, bluntly stated, ``I wouldn’t have written it that way, but it will do,'' and immediately pivoted right back to our latest research with his enthusiastic: ``But you {\em have} to see this...''

He loved to travel to learn new mathematics, and since he needed help with logistics, I was a frequent companion on his foreign trips. He was, of course, legendary for beach meetings. One day, sitting on the sand at La Jolla Shores with Michele D’Adderio and myself, he decided that the only thing better than doing math on the beach in La Jolla was doing it in Hawaii. Upon returning home, he negotiated a deal with his wife: he could take the two of us to Hawaii if the family got a dog. Years later, he would still get a goofy grin describing the trade, admitting his wife got the better end of the bargain: the trip was long over, but they still had the dog. Still, I always thought he was wrong, since Michele and I came out the best. 

Adriano had a remarkably diverse group of students for his era; roughly a quarter of his PhDs were women. He was somewhat less progressive on work-life balance, looking quite sad to meet my then-boyfriend (now husband) when I graduated. He brightened up considerably when I explained that because my boyfriend was still finishing his degree, I would have an excuse to return regularly to San Diego to do math. While gregarious about his own kids, he also frequently bemoaned how difficult it was to balance motherhood and a career, though usually in a manner that suggested he simply hoped mathematics would be my sole child.

If there is one trait that truly defined Adriano, it was how fiercely he fought for his students. Many of my siblings skirted graduate school rules because Adriano interceded; they are now top researchers in our field. He frequently gave up his own NSF summer salary to support our research and travel. When we traveled together, he would loudly boast about my ``amazing'' achievements to international audiences long before mentioning his own work. He was a ferocious champion on the job market, giving us a rare brand of confidence and aggressively opening doors that gave our careers a profound head start. May his {\em joie de vivre}, sharp mathematical insights, and endless passion for supporting young mathematicians live on in the many students he championed!

\section*{Luc Lapointe}
I first met Adriano in 1995 while I was doing my PhD in Physics at the Université de Montréal.
I had never seen anyone before give a talk with such passion and enthusiasm, which was absolutely inspiring.
And he was not afraid to make bold statements! I still remember that he said something along the lines of ``The analysts are the bad guys, who cares about convergence?'', which truly resonated with me as a physicist.

Adriano was always extremely generous and enthusiastic about the work of young researchers, and his encouragement about my work was a huge boost to my morale. One time I was in a new office and Adriano was looking for me.  I can still hear his voice asking ``O\`u est Lapointe? O\`u est Lapointe?''.  Nothing can be better than having someone of his stature making you feel like you are special  (as many have shared during his online tribute, he had a remarkable ability to make you believe in yourself).

\begin{wrapfigure}{l}{0.5\textwidth}
    \centering
    \includegraphics[width=.95\linewidth]{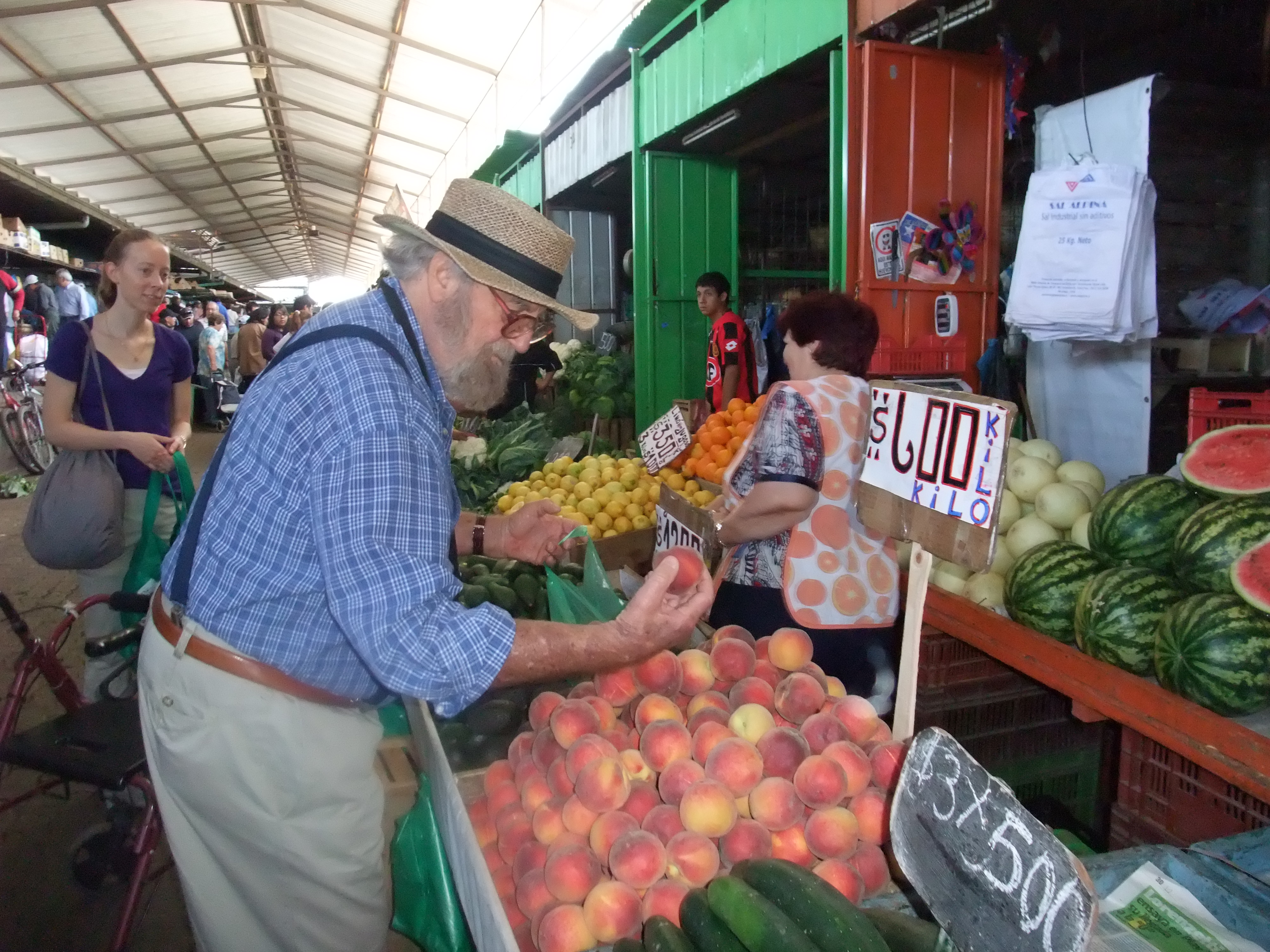}
{\footnotesize Adriano inspects fresh fruit in a Chilean market}
\end{wrapfigure}

But Adriano's generosity went far beyond just words of encouragement. The following year, he invited me to UCSD for I.G. Macdonald's lectures on double affine Hecke algebras. I stayed at his home, where I had the chance to meet his wife, Diane, who has always been incredibly kind to me, and his lovely daughter, Gabriela.  There was instant familiarity with Adriano.  He installed me in his office, where we worked in the mornings. And then we spent the afternoons working at the beach.  At some point, I even went with him to pick up I.G. Macdonald at the airport!  

Adriano's influence on my career and research was profound. I spent the second year of my postdoc at UCSD (he convinced me to first spend one year in Marne-la-Vallée). At that time, there was great excitement and development about Macdonald polynomials and  the $n!$ conjecture. I remember being especially thrilled about his Science-Fiction paper with Fran\c{c}ois Bergeron, which posited that the intersection of $k$ Garsia-Haiman modules had dimension $n!/k$ (the idea of understanding Macdonald polynomials through refinements has motivated me ever since). While I never formally collaborated with Adriano, I learned so much from him, especially in the art of ``manipulatorics'' that he mastered like no one else. 

I also always felt his genuine excitement for my work, and  his support often went beyond simple encouragement (he dedicated for instance countless hours shaping our first paper with Alain Lascoux and Jennifer Morse on atoms, contributing invaluable insights to the introduction and overall presentation). My time at UCSD also meant a lot of time working at the beach (which I never fully enjoyed) and a lot of dinners at his house (which I always greatly enjoyed).

I saw him less often after I moved to Chile in 2002.  He visited Talca only once (he was very pleased that the scenery reminded him of his youth in Tunisia, but above all he loved the avocados!). I remember a poignant moment during a 2006 visit to San Diego. He sadly confided that he was afraid that the end was coming, that he wouldn't be able to do mathematics forever as he wished. I was extremely moved, yet I couldn't find any words of comfort. Forever was probably a lot to ask, but never did I imagine that he would keep doing first-rate mathematics for more than 15 years.  It seems that his work  even blossomed from that point on, and he belonged to a growing community working on fascinating problems in $(q,t)$-combinatorics such as the shuffle and Delta conjectures. Even more amazingly, he had 7 or 8 PhD students in a row getting NSF postdoctoral grants while he was in his 80's and 90's.

I have been fortunate to have met some incredible people, but Adriano was in a league of his own. Everything about him was larger than life. His energy, his enthusiasm, his passion, his humor, his cooking. I especially enjoyed my conversations with him as he had such a unique perspective  on both life and  mathematics (as well as some of the best stories I have ever heard!). On a deeper level, he was the mathematician that I most wanted to impress, and I will terribly miss this driving force.

\section*{David Little}

Like most of Adriano's former students, I have so many fond memories of being a student in his graduate Applied Algebra and Combinatorics courses and as his teaching assistant for his undergraduate courses in Cryptography, Computer Graphics, and Animated Statistics.  And then there's all the mathematics lessons at the beach, great food and conversations at his house, and traveling to conferences with him, both near and far.  But instead of focusing on any one specific memory, I'd like to share a couple of general observations I've had about Adriano. 

First, I never saw him to be a competitive person.  He shared so many of his ideas and conjectures with his students, and I never once heard him express any regret for doing so when someone else solved a problem he had worked on.  As far as I could tell, he just wanted to see beautiful mathematics and he didn’t care if it came from him or one of his students.  Of course, most of the time the beauty was of his making, but in those cases when it wasn't, he would usually find a way to make your work better. I recall numerous conversations with Adriano where he told me about his latest graduate students and the amazing work they were doing.  He was so proud of all of them.

Second, he genuinely cared for his students. He invited us into his home and made us feel like part of his family. Over the years, through numerous conversations with colleagues who didn't have as much contact with their advisors, I've come to appreciate how rare that is. When I was on a plane to San Diego for a conference in honor of Adriano, I just happened to be seated next to another mathematician.  He wasn't going to Adriano's conference, but when I told him about it, he was truly amazed at how many of his students had kept in touch with him through the years and were willing to travel from all over to see him again.  I hope Adriano knew how much his students cared about him. 

Adriano was extremely generous, creative, supportive, passionate, and not just about mathematics, but about living life to its fullest. He was a one-of-a-kind mathematician and person.  He built an incredible community of mathematicians, and I consider myself extremely lucky to be a small part of it.  I will always be indebted to Adriano for the wonderful world of mathematics he shared with me and for the undivided attention he gave me whenever I had something to share in return.  He will always be missed and never forgotten.

\section*{Arun Ram}
(1)  A reason to go to UC San Diego for a PhD:\\
That man had fire,\\
That man had culture\\
That man had passion\\
That man sparked electric energy\\
And I could learn mathematics from that man.\\
And that, after all, was the point.
\vskip .1in

(2) In graduate school, I failed my quals for the nth time.  
The rules said I had to be kicked out out of the PhD program.
Garsia rose like a great Durga,
riding a tiger and wielding weapons with many arms,
fighting for my cause.  
It didn’t make him many friends among his colleagues, but he won.
I didn’t get kicked out of graduate school and
my mathematics career had another chance.
I was immensely immensely grateful,
but never dared to bring it up
for fear that Durga might rise again.
\vskip .1in

(3) One Sunday I came up the elevator in the math building.
The elevator door opened and  I could hear yelling down the hall.  
In trepidation and curiosity I edged my way hesitantly toward the commotion.  
As I got closer I realised that it was Adriano (a world expert in Schur functions)
in Jeff Remmel’s office (another world expert in Schur functions)
screaming ``What’s the definition of a Schur function!!”  
I slunk cautiously and quietly to the stairwell and went quickly down the stairs.
I learned an important mathematical lesson that day:
On Sundays and Wednesdays,
the definition of a Schur function might be different.
It has been a very useful lesson in my work since that time.
\vskip .1in

(4) Garsia was always happiest when you were engaged and working
together with him on what he was working on.
I listened avidly and intensely but I worked on other things.  
One of my greatest successes in mathematics was
when I  managed to get Garsia to work,
for about 20 hours, on what I was working on.  
I had given a talk in the seminar, and he got interested
and spent the evening and the next morning
working through what I had presented,
amplifying and processing it from his own point of view.    
The next afternoon he gave me two hour lecture explaining
how it could all be ``derived from scratch''.  
He ended the session with an apology
that he couldn’t possibly be pulled away any longer
from his really important work on kicking boxes
(it was the time of Garsia-Procesi modules).  
And then he started another long lecture on kicking,
which got interrupted, because
the sun was getting low,
the beach was getting cold,
and it was time for Garsia to go home and cook dinner.
\vskip .1in

(5) Wonderful memories of cooking and dinners
outside on the back patio of Garsia’s house.  
While barbecuing, Adriano would tell stories of the farm in Tunisia --
cooking over the open fire,
finding unexploded bombs left over from  WWII around the property,
and squeezing the roasted eggplants and peppers
so that the soft juicy interior  oozed out of the skin and
mixed with the  roasted garlic that had been extracted from the coals.  
We would eat well and  thoroughly
(even for a young graduate student),
while listening to stories of life and
trying not to get distracted by the beautiful sunset
and evening twilight, looking over the bay,
the night lights of the city of San Diego,
the ocean and the shining moon.
Since that time, Tunisian roasted vegetable sauce remains
a staple in my family’s household -- over pasta, over sandwiches, over anything.
\vskip .1in

(6) The phone rang at 7am on Saturday morning.  ``Arun, are you up?''
(Like any proper graduate student
I had been planning to sleep until 1:30pm after a late night out.).
Before I could muster an answer, the voice continued:
``I’ve been waiting for hours until it was a suitable time to call,  
you must come over right away, it’s fantastic, I have to tell you all about it''.
I knew that meant that there was a wonderful breakfast in store for me
and a couple of hours of exciting mathematics explanations.  
``I need to take a shower, but then I’ll be right over'',
and sure enough 45 min later,
I was enjoying an aromatic double espresso, croissants and eggs and fruit and cheese
and listening to the expostulations of how to modify Macdonald polynomials
into the most fantastic object ever observed by humans.
\vskip .1in

(7) 30 years later Adriano wrote to me a typical message about a paper I had recently sent him:
\begin{Verbatim}[fontsize=\small]
I have tried to read the papers you sent me
unfortunately they were not written for me...
I need more definitions,
I could not even figure out what Macdonald polynomials you are involved with..

I only work with the symmetric group case and
even in that case
I work only with the modified ones...

I have enough beautiful open problems in that case to work with
to start getting involved with conjectures in more general cases.

Still I wanted to see to what extent your world was connected to mine..

Can you tell me at least what Macdonald polynomials your paper is concerned with?
\end{Verbatim}
\vskip .2in

I responded 6 months later:
\begin{Verbatim}[fontsize=\small]
Dear Adriano,
So, in your last email to me you gave me a problem.
I am very slow and it took me some time to work it out in a simple way.
I’m still old fashioned and I wrote a handwritten letter.
The scan is attached below
(I’ve also sent the original to your home address.)

I’ll be interested to hear if you have some thoughts.
\end{Verbatim}
\vskip .2in

His email response came the next day:
\begin{Verbatim}[fontsize=\small]
Dear Arun,

Thank you for your great email!!!

I have read your very impressive hand written letter in great detail.

I am specially amazed by your hand calculations at the end of your letter...

Do not short change yourself, by saying that you are slow to get things done...
Speed is not essential, specially if it is not supported by persistence.
There is a famed fable by Lafontaine "le lievre et la tortue."
The tortue wins the race, and not the lievre, since she is more persistent.

I wish very much I could get you interested in working
in some of the open problems in my area.
You are getting very close indeed,
closer than you have ever been…
but not quite yet.
\end{Verbatim}

Garsia was always happiest when you were engaged and working together with him on what he was working on.
It was, in so many ways, the most fantastic, amazing and thrilling way to tantalize your mind.

\section*{Marino Romero}

I have tried many times to write something that encapsulates how Adriano has affected my life, and each time I feel I either am missing something or I am talking too much about myself. In many ways, my story is common around Adriano--he has helped so many people and been an important part of other people's successes. But I will do my best to explain why Adriano and his family are so dear to me. 

When I was an undergraduate student at Cal Poly SLO, my undergrad advisor, Tony Mendes, showed us what I would now call ``UCSD combinatorics''. Among the topics we learned, we saw the Garsia--Milne Involution Principle and Remmel's Bijection Machine. That was one of the first times I ever felt the magic of mathematics. I guess it may be why I always try solving problems with sign-reversing involutions--I just found the method so beautiful.

It was here that I first learned, through stories, about Adriano's ability to keep talking about his current mathematical interests, even when the recipient was tired of listening; and here I also learned about how he loved to work on mathematics at the beach. It became like a myth to us, something I never imagined witnessing.

When I was graduating, I knew I wanted to continue studying mathematics, but I did not get into any graduate program. Most likely my statement was horrible, I did not know I needed to prepare for the GRE, and there are probably many other reasons. But I had planned to continue studying by taking a class at UCSD while I worked in whatever job I could get my hands on. Tony Mendes told me Adriano would be teaching the ``Applied Algebra'' course, so we thought that would be a great class for me.

By the following year, Adriano had become my advisor and I had passed my first qualifying exam (all thanks to both Adriano and Jeff Remmel). But most importantly, I got to go on my first trip to the beach with Adriano. One of the highlights was that we immediately got into a car accident when we got there, when a big truck hauling boards and kayaks backed up into my passenger side. He quickly exchanged information, waved it off, and moved right into working on some mathematics with me. Ever since then, I would go work on mathematics at the beach by myself.

Many events were like this with us. When he first took me to Montreal, he got a big wound on his leg during our connecting flight. This led to a sequence of events that I always remember and laugh about. Besides having to assist him, when we got to his condo, he had me carrying a washer and dryer around because insurance people had to come and check recent flooding damage. I don't know exactly how many times I had to move them back and forth, but luckily I was a lot more physically fit back then. 
The following day, he left because his wound was getting infected, but I continued there, in the coldest weather I had ever experienced in my life. 

There was a time where we went to a store just to get the bread that he liked. His driving always involved a little ``pep,'' as I remember him describing cars as not having enough ``pep.'' All of the sudden a car tries to cross, and Adriano pulls a wonderful maneuver to avoid an accident. With his age, I had a hard time believing how well he handled it. In general, however, there were a number of times that I wanted to drive instead.

Many of the times I met with Adriano, the day felt like a quest. For example, if we wanted to make the perfect pasta sauce with meat, we had to travel to find the perfect steaks. And after the voyage, we would return to what was important to us--solving a mathematical problem. 

I felt a connection that is difficult to explain. Sometimes, we would spend a long time yelling at each other, me trying to convince him that what I was saying was correct. In the end, he would agree then say something like, ``Ahh... Want to go to Rubio's (for some fish tacos)?'' Several times he would say to me, ``you know us Latins, we just get very emotional.'' These kinds of controversial statements may sound bad, but he was right in many ways. We argued because we cared about the subject, and I truly appreciated the effort he took to hear what I had to say. 

I am also very grateful of his family. Gabriella, Diane, Adrian, Maya, all who were incredibly welcoming. They made me feel at home while eating dinner with them. These are the memories I hold most dear in my life. Sure, I learned a lot of mathematics, beyond anything I dreamed of. But I also felt a comfort while doing so that I do not think I would have experienced with another advisor.

\section*{Emily Sergel}
I met Adriano at the UCSD open house in 2011. From this one meeting, I knew I wanted him to be my PhD advisor. I am so grateful for the time I got to spend with him, and I hope the stories below can show a small fragment of what made him such an amazing person.

According to Adriano, doing math was the solution to all of life's problems, and he couldn't imagine that anyone would feel differently if they gave math a try. When a friend announced he and his wife were getting divorced, Adriano decided to teach the wife some math. He said that doing math would make her happy and then they wouldn't need to get a divorce.

He truly believed that anyone could succeed at math with the right attitude. He often told a story about redoing his driveway. Rather than pay a crew to demolish the old one, he threw a party for faculty and students. He invited everyone to have a turn using the jackhammer. Many would try for a minute or two and then give up. But some stuck with it and kept hammering, eventually breaking through the concrete. The students who stuck with it, Adriano claimed, were the ones who later wrote the best PhD theses.

He was optimistic, proactive, and resourceful, tackling anything that caught his interest with determination. He taught himself how to cook without the guidance of recipes, based only on the memory of how his mother’s food had tasted in his childhood. He learned programming while teaching it, back when the UCSD Computer Science department was new and in need of instructors. He said many times that teaching something was actually the best way to learn it yourself.

During my first year, Adriano's bad back prevented him from standing for a whole class. So, he invited me to go up to the board instead. He told me what to write and asked me questions. We filled the class with a mix of lecture-by-proxy and Socratic dialogues. We ended up doing this together for all his graduate courses, and he continued the same way with his next student. He was right – this was an excellent way for me to learn, and it gave me a unique perspective on how to teach mathematics too.

Adriano is one of my favorite speakers. His math talks were chaotic, passionate, and inspiring. When he made slides, he always included colorful diagrams and animations. Even when talking about something for the fiftieth time, he would still get visibly excited, practically jumping up and down on the punchlines. His enthusiasm was irresistible and drew countless others to work on the problems that he loved.

When I started graduate school, Adriano was already in his 80s. Some people thought I was crazy to choose such an old advisor. But the people who said this had never met Adriano. Even with a bad back and a walker, anyone who knew him couldn't really think of him as old. The idea that he may retire soon was ridiculous. He was full of life, and full of passion for math and food and music and travel. When he did eventually retire – many years after I graduated – he kept doing mathematics, working with students, swimming in the ocean, cooking delicious meals, and living his life with joy.

\section*{Michelle Wachs}

\begin{wrapfigure}{r}{0.4\textwidth}
    \centering
    \includegraphics[width=0.95\linewidth]{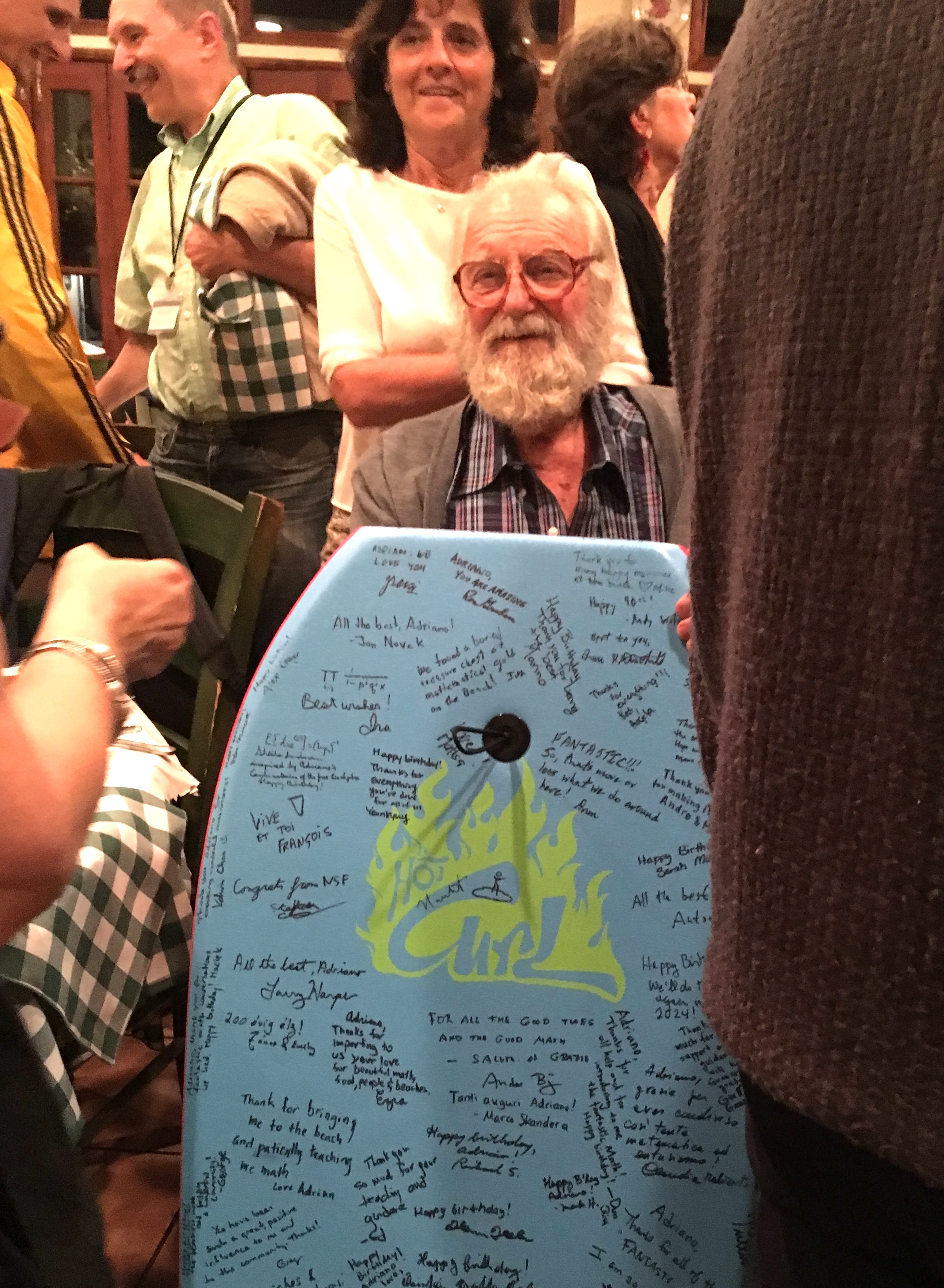} 
 {\footnotesize Adriano and the Boogie Board at Garsiafest (with Luisa Carini
 and Andre Broido also pictured), Scripps, 2019}
    \end{wrapfigure}

I was one of Adriano’s first combinatorics students at UCSD in the mid-1970s. I met him during my second year of graduate school when I took his Real Analysis course. He was an extraordinary teacher, and since he was also renowned as a leading researcher in analysis, I hoped he might become my advisor in that field.  Before I even got up the courage to ask him, he invited me to his office to share  ideas he had about simplifying the proof of the 
Hu--Tucker algorithm for optimal binary search trees. I was surprised: this wasn’t real analysis! In fact, it was combinatorics, though I didn’t know what that was at the time. I also didn’t know that Adriano was transitioning from analysis to combinatorics, where he would become one of the world's foremost leaders in algebraic combinatorics for almost 50 years.

I worked on Adriano's binary search tree problem and had a breakthrough:  his ideas led, not to a simpler proof of the Hu--Tucker algorithm, but  to a  new simpler algorithm for constructing optimal binary search trees. 
I was so excited to show him what I had done, I went to the university and sat on the floor outside his office door for what felt like hours waiting  for him to arrive. (This was long before email, cell phones, or texting.)  Little did I know, he was at his favorite place -- the beach -- which also happened to be my favorite place -- what luck!   

The beach was Adriano's second office. 
His students, collaborators, and many distinguished visitors would find him at his  encampment on the beach (first in Del Mar, later at La Jolla Shores) where he would share his beautiful mathematical ideas and listen to theirs. He also shared his passion for riding the waves.    He got a Morey Boogie Board when they first came out in 1975, and after seeing his, I went right to the store to get one too. In the years to come he would bring a stack of them to the beach for his visitors to use.  At his 91st birthday conference — held appropriately at the Scripps Institute of Oceanography overlooking the beach — a boogie board signed by all the participants was presented to him.  The NSF program director remarked in his speech at the conference that  Adriano, at $91$, was the oldest active principal investigator on an NSF grant in any scientific field. 

\begin{wrapfigure}{l}{0.5\textwidth}
    \centering
    \includegraphics[width=0.95\linewidth]{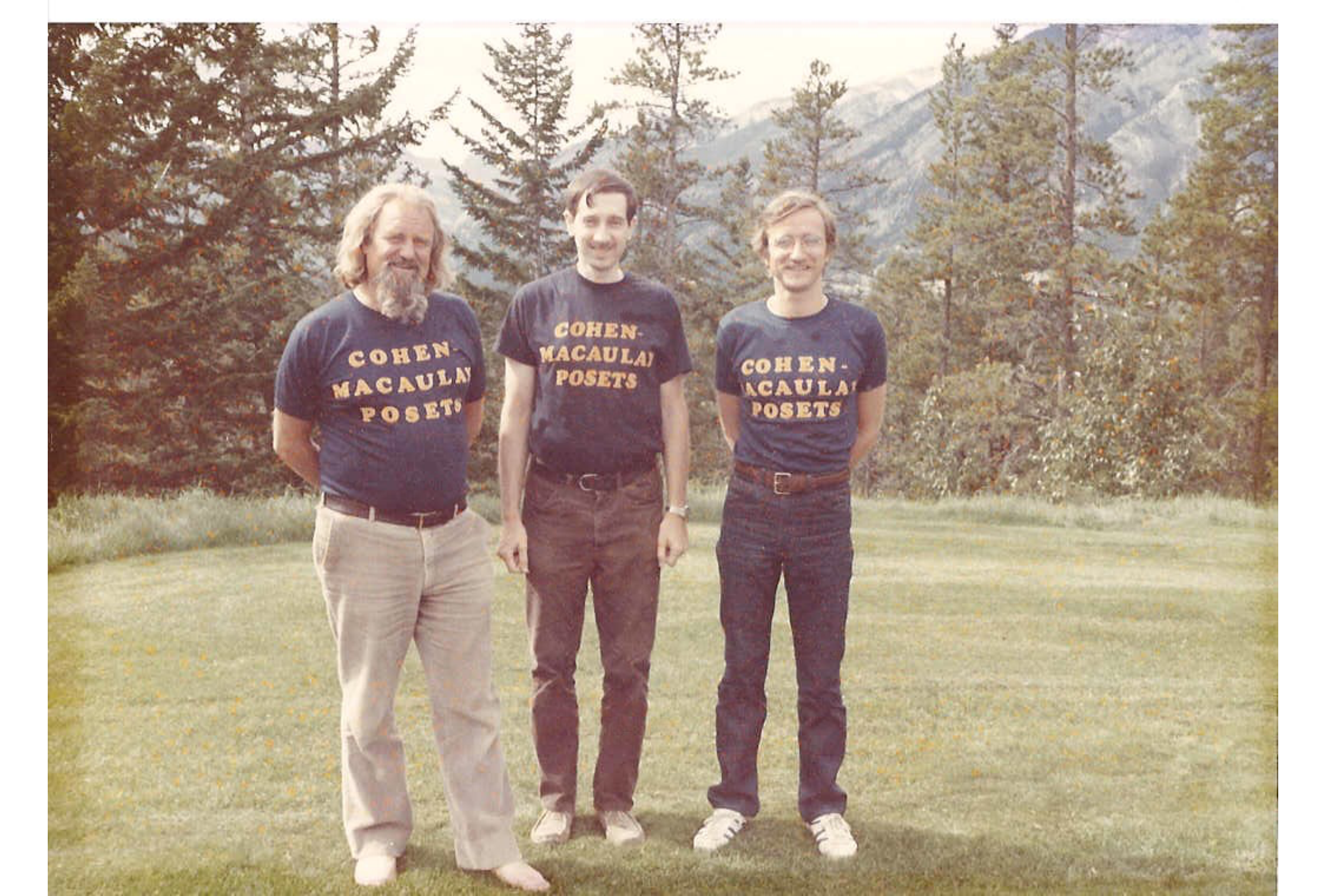} 
  {\footnotesize Adriano, Richard, Anders at Banff Conference on Ordered Sets, 1981}
\end{wrapfigure}

 From the mid-1970s to the late-1980s, Adriano's pioneering work helped shape the field of algebraic combinatorics and made him  one of the field's true founding leaders. So many of the ideas he introduced during that  period continue to influence the field today. I’d like to share just a few threads from that time that have had a deep and lasting impact on my own work.
 Adriano’s early research on 
$q$-analogs and symmetric functions — much of it done in collaboration with Ira Gessel and Jeff Remmel — laid the groundwork for his later celebrated work with Mark Haiman on Macdonald polynomials and diagonal harmonics. Those early papers -- as well as his subsequent papers --  were also a key source of inspiration for my own research on 
$q$-analogs and symmetric functions throughout my career. His work on Cohen-Macaulay posets, including his influential survey chapter  with Anders Bj\"orner and Richard Stanley, played a central role in the early development of topological combinatorics,  a field I’ve been fortunate to be part of since  early in my career. And Adriano’s beautiful work on the Lie representation of the symmetric group continues to resonate with me. It sparked ideas that led to some of my own early research and has remained a source of insight and motivation in my work on this topic ever since.

Adriano was an exceptional advisor -- extremely generous with his time and ideas. He had a rare gift for recognizing and nurturing the individual talents and abilities of his students.   Many of us feel we owe what ever success we have achieved to his mentorship.  In addition to sharing his great ideas, he was supportive and passionate.  His excitement about our work was a great confidence booster.  He and I often had spirited debates when discussing math,  but we always resolved things amicably.   The good thing was that I always felt free to express myself openly and honestly. 

\begin{wrapfigure}{r}{0.4\textwidth}
    \includegraphics[width=0.95\linewidth]{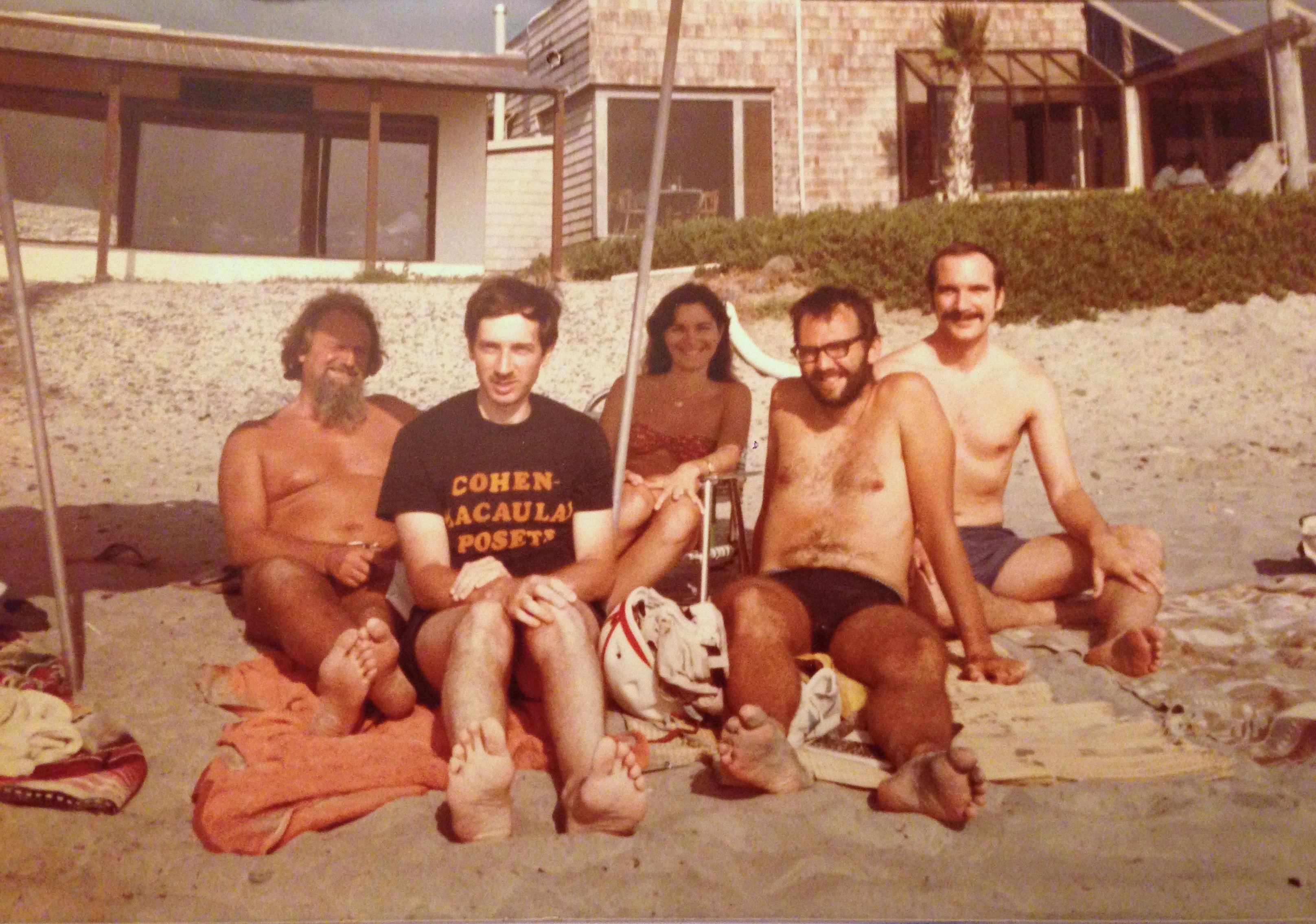}\\
  {\footnotesize UCSD Combinatorics Seminar, Del Mar, 1981,
Adriano, Richard, me, Dennis, Jeff}
\end{wrapfigure}

Adriano was there for his students long after they graduated.    For many years after completing my Ph.D., my husband (a fellow UCSD Ph.D. -- in differential geometry)  and I would spend the summers visiting UCSD, and we even spent a full academic year there as visiting associate professors.   
Adriano       helped make this possible by being a wonderful mathematical host and finding  places for us to stay.   
During those visits, there was a vibrant atmosphere in algebraic combinatorics, with seminars given
by distinguished visitors such as Richard Stanley, Anders Bj\"orner, Alain Lascoux, and Dennis Stanton, to name a few. 
During my first summer visit, Adriano suggested to Anders that he discuss with me his problem
on shellability of Bruhat order.   This led to a lifelong collaboration and friendship with Anders
for which I remain deeply grateful to Adriano. This was typical Adriano—he fostered many enduring mathematical relationships.

My husband Greg and I took some memorable trips with Adriano and his wife Diane, both mathematical and personal.  We had a great time in Rome with Adriano teaching us all how to behave in Italy and how to sniff out good restaurants. Our young children, Brian and Gabriela, were inseparable during our stays in 1992 at Mittag-Leffler and in 1994 at Garsiafest in Taormina. Whether at home or away,  Adriano loved to entertain his colleagues and friends by treating us to unforgettable Italian and Tunisian meals, which he cooked himself.

To so many of us, Adriano was more than a mentor, colleague, or friend—he was family.  Once, during a short visit to San Diego, he and Diane invited me to stay with them. At the time, their two-bedroom country home  was full, with Diane’s brother occupying the second room.  So my ``room'' was the front porch.   I left my suitcase in the living room while I slept out on the porch.  My suitcase was a  beat-up leather suitcase with a handle that had become partially detached.  When I came in the next morning, to my surprise, Adriano  was busy sewing the handle back on with a giant needle and thick thread. I was so touched. What could be sweeter and more fatherly than that?  

\begin{center}
 \includegraphics[width=0.7\linewidth]{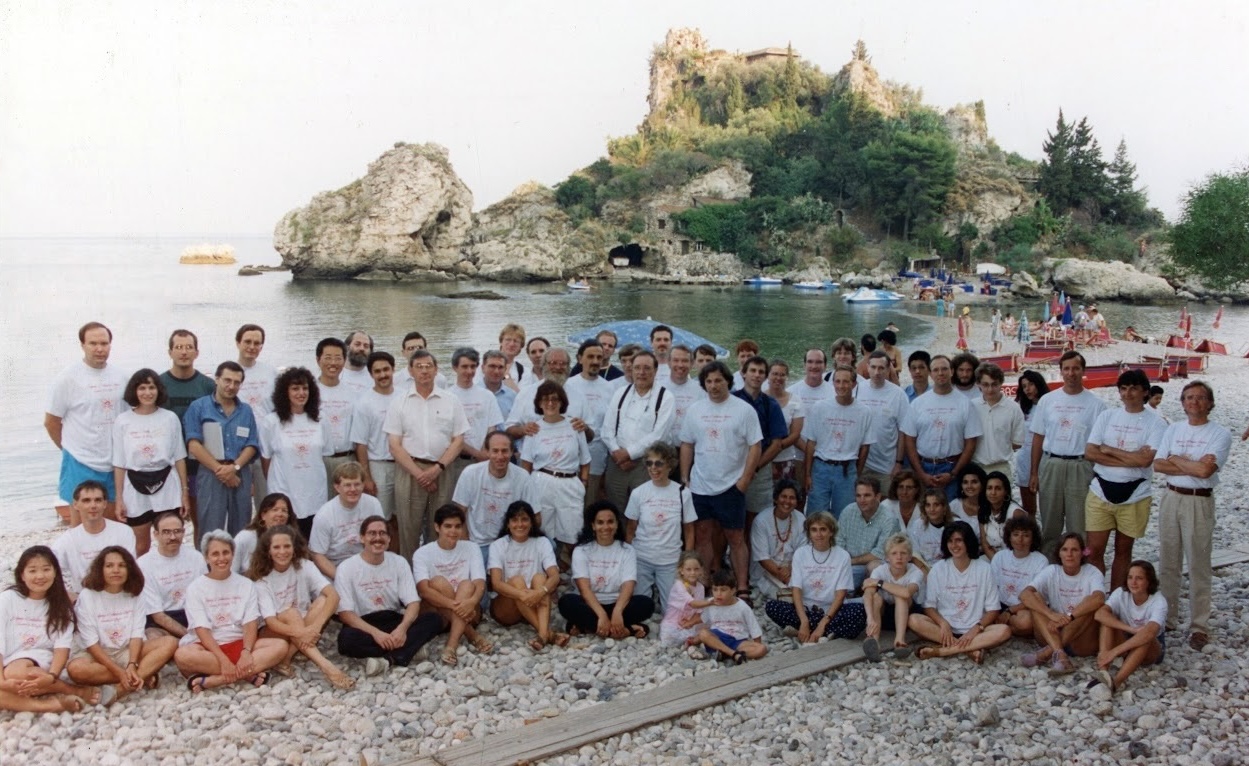} 
 \\ {\footnotesize Garsiafest, Taormina, 1994}
\end{center}

Adriano’s passion for mathematics—and for life—was legendary.  His work continues to have a strong influence on my work and that of so many others.  What a profound privilege it has been to have known this remarkable man for so many years.

\section*{Guoce Xin}

Adriano Garsia was more than a mathematician—he was someone who truly loved mathematics for what it was. That passion stayed with him all his life. He worked on the Shuffle Conjecture for years, and every step forward made him genuinely happy. When Anton and Mellit finally solved it, he wasn't disappointed at all; instead, he was eager to go through their proof and understand it himself.

In 2006, he found me online, and we began working together on the sdd5 problem. His sharp logic and energy for math were so vivid that in those days I never thought about his age. It was only when I visited him in 2014 that I saw him help my family settle down, just like a grandfather—and by then I did know his age, but he never seemed old to me. He always encouraged me and his students, even for small progresses. He liked to write things clearly himself, but his only demand from us was that we explain ourselves clearly. I learned a great deal from him, and I have long considered him my second advisor, alongside my PhD advisor Ira Gessel. One piece of advice he gave me still stays with me: when you tackle a hard problem, don't give up too soon—but if you're stuck for two weeks, it's fine to put it aside.

He loved the sea. Sometimes we would talk about math by the shore, looking out at the wide horizon that seemed to stretch on forever, much like mathematics itself.

\section*{Mike Zabrocki}

One of the first upper division courses that I took
as an undergraduate student in the very early 90's
was a course on cryptography from Adriano Garsia.
I learned years later that Adriano taught the course because
he hated teaching large first year calculus classes and so in order to get out
of them he created several courses that allowed him to satisfy the department that
he was teaching to undergraduate students while allowing him to be more creative in the
subject matter that he taught.  Not all of the courses that he proposed were popular enough
to attract a large enough enrollment survive, but this one was always full.

\begin{wrapfigure}{r}{0.4\textwidth}
    \includegraphics[width=0.95\linewidth]{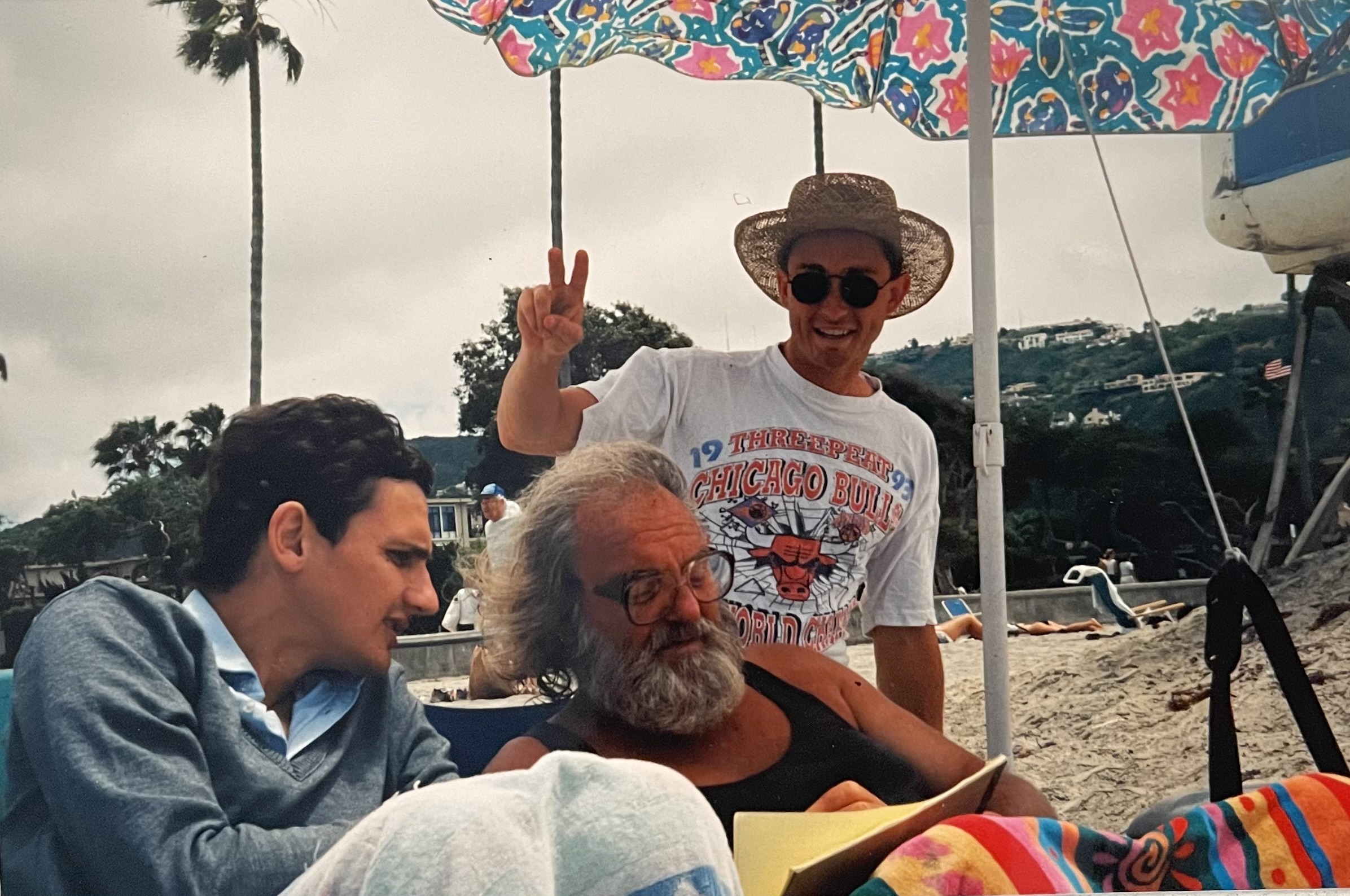}\\
  {\footnotesize Christian Krattenthaler visiting Adriano at the beach while
  we were goofing off.}
\end{wrapfigure}

When he introduced the Vigen\`ere cipher as one of the early topics in the course,
I wrote a program in Pascal that would guess at the encryption key
given the cyphertext alone using single letter statistics.
The program wasn't revolutionary, but I was quite
proud of myself since it wasn't something that he had showed us how to do in class.
I sent an email to him showing him what I had done and
this triggered an interaction that I would later in life recognize as typical of
Adriano when he got excited about mathematics:
a burst of enthusiasm and praise and showering of attention.
The next class, in his own unique style, he showed off what I had done but in a way that
left me baffled.  He had translated my program into 5 lines of APL
code that only he could possibly read and then showed it to the class and attributed
the program to me.

The following summer he hired
me as a research assistant.  I remember him driving me down to the beach in his van where
we met with Mark Haiman who showed me how projection operators could be used to reduce
computations to calculate the
dimensions of linear spans of derivatives of a polynomial in two sets of variables.
His patience with me was emense because I had just finished my second
year of courses and didn't have the slightest understanding of what a vector space was,
let alone why we were trying to compute the dimension of them.

I was already trapped in the gravity of Adriano's personality when I started graduate school,
but it would be a few years before I was aware of that.  I was his teaching assistant for his
cryptography class and several other courses that he taught
and I had little chance of reaching escape velocity once I took his course in applied algebra.

Adriano's perihelion was at the beach.   In the days before cell phones,
those of us in his orbit would call him
in the morning while he was at home and ask what time he expected to be there.
On any given afternoon, a former students or co-author who happened to be
in town to know that he would be there and to drop by to say hi and talk math.

When I wrote my thesis I included the following paragraph in the acknowledgement
section that still captures my feelings about my experiences with Adriano.

\begin{quote}
I feel that I have been changed the most by the influence of my advisor, Adriano.
I've loved all the places I have traveled, the people that I have met, the things I have
learned to cook. These are things that I will take away from graduate school that mean
much more to me than the mathematics that I've learned from him. The `where,' `how,'
and `what' mathematics I do is definitely dependent on his influence and is visible in
this document. It is less visible that he has changed my personality. Many of the things
that I most like about myself I have learned from him. `I am the student of my advisor.'
\end{quote}

What still inspires me about Adriano is how the people that did mathematics with him
became part of his family.  He integrated us all into his life in a way that
you can sense from the anecdotes presented by his colleagues and former students.
\end{document}